\documentclass[a4paper,12pt,reqno]{amsart}

\usepackage{amsmath,amsfonts,amsthm,amssymb,color}
\usepackage[T1]{fontenc}

\def\DD{{\mathcal D}}
\def\VV{{\mathcal V}}
\def\UU{{\mathcal U}}
\def\EE{{\mathcal E}}
\def\HH{{\mathcal H}}
\def\TT{{\mathcal T}}
\def\RR{{\mathcal R}}
\newcommand{\bt}{{\bf T}}
\newcommand{\btt}{{\bf t}}
\newcommand{\bss}{{\bf s}}

\newcommand{\der}{\delta}

\newcommand{\hro}{\hat \rho}

\newcommand{\hz}{\hat z}

\newcommand{\id}{\mbox{Id}}
\newcommand{\iou}{\int_{0}^{1}}
\newcommand{\iot}{\int_{0}^{t}}

\newcommand{\ist}{\int_{s}^{t}}
\newcommand{\norm}[1]{\lVert #1\rVert}

\newcommand{\ott}{[0,T]}
\newcommand{\ou}{[0,1]}

\newcommand{\1}{{\bf 1}}

\newcommand{\D}{\mathbb D}
\newcommand{\R}{\mathbb R}
\newcommand{\N}{\mathbb N}

\newcommand{\cb}{\mathcal B}
\newcommand{\cac}{\mathcal C}

\newcommand{\cd}{\mathcal D}

\newcommand{\ch}{\mathcal H}

\newcommand{\cj}{\mathcal J}

\newcommand{\cu}{\mathcal U}

\newcommand{\cz}{\mathcal Z}

\newcommand{\al}{\alpha}
\newcommand{\ep}{\varepsilon}

\newcommand{\ga}{\gamma}
\newcommand{\gaa}{\Gamma}
\newcommand{\ka}{\kappa}
\newcommand{\la}{\lambda}
\newcommand{\laa}{\Lambda}

\newcommand{\si}{\sigma}

\newcommand{\vp}{\varphi}

\newcommand{\lp}{\left(}
\newcommand{\rp}{\right)}
\newcommand{\lc}{\left[}
\newcommand{\rc}{\right]}
\newcommand{\lcl}{\left\{}
\newcommand{\rcl}{\right\}}
\newcommand{\lln}{\left|}
\newcommand{\rrn}{\right|}
\newcommand{\lla}{\left\langle}
\newcommand{\rra}{\right\rangle}

\newcommand{\bean}{\begin{eqnarray*}}
\newcommand{\eean}{\end{eqnarray*}}
\newcommand{\ben}{\begin{enumerate}}
\newcommand{\een}{\end{enumerate}}
\newcommand{\beq}{\begin{equation}}
\newcommand{\eeq}{\end{equation}}

\newtheorem{theorem}{Theorem}[section]

\newtheorem{corollary}[theorem]{Corollary}

\newtheorem{lemma}[theorem]{Lemma}

\newtheorem{proposition}[theorem]{Proposition}

\newtheorem{hypothesis}{Hypothesis}

\theoremstyle{remark}
\newtheorem{remark}[theorem]{Remark}
\newtheorem{remarks}[theorem]{Remarks}

\begin{document}

\title{Malliavin calculus for fractional delay equations}

\author{Jorge A. Le\'{o}n \and Samy Tindel}

\address{
{\it Jorge A. Le\'{o}n:}
{\rm Depto. de Control Autom\'{a}tico, CINVESTAV-IPN, Apartado Postal 14-740, 07000 M\'{e}xico, D.F., Mexico}.
{\it Email: }{\tt jleon@ctrl.cinvestav.mx}
\newline
$\mbox{ }$\hspace{0.1cm}
{\it Samy Tindel:}
{\rm Institut Élie Cartan Nancy, B.P. 239,
54506 Vandoeœuvre-lès-Nancy Cedex, France}.
{\it Email: }{\tt tindel@iecn.u-nancy.fr}
}

\keywords{Delay equation, Young integration, fractional Brownian motion, Malliavin calculus.}

\subjclass[2000]{60H10, 60H05, 60H07}

\date{\today}

\thanks{J.A. León is partially supported by the CONACyT grant 98998. S. Tindel is partially supported by the ANR grant ECRU}

\begin{abstract}
In this paper we study the existence of a unique solution to a general class of Young
delay
differential equations driven by a H\"older continuous function
with parameter greater that $1/2$ via the Young integration setting. Then some estimates of the
solution are obtained, which allow to show that the solution 
of a delay differential equation driven by a fractional Brownian motion
(fBm) with Hurst parameter $H>1/2$ has a $C^\infty$-density. To this purpose, we use 
Malliavin calculus based on the Fr\'echet differentiability
in the directions of the reproducing kernel Hilbert space
associated with fBm.
\end{abstract}

\maketitle

\section{Introduction}

The recent progresses in the analysis of differential equations driven by a
 fractional
Brownian motion, using either the complete formalism of the rough path analysis
\cite{CQ,Gu,LyonsBook}, or the simpler Young integration setting \cite{NR,ZA},
 allow to study some of the basic
properties of the processes defined as solutions to rough or fractional 
equations. This global program
has already been started as far as moments estimates \cite{HN}, large
 deviations \cite{LQZ},
or properties of the law \cite{BC,NNRT} are concerned. It is also natural to consider some of the natural generalizations of diffusion processes, arising
in physical applications, and see if these equations have a counterpart in the
fractional Brownian setting. Some partial developments in this direction concern pathwise type PDEs, such as heat \cite{DGT,GLT,GT,Te}, wave \cite{QT} or Navier-Stokes \cite{CQT} equations, as well as Volterra type systems \cite{DT1,DT2}. As we shall see, the current paper is part of this second kind of project, 
and we shall
deal with stochastic delay equations driven by a fractional Brownian motion
with Hurst parameter $H>1/2$.

\vspace{0.2cm}

Indeed, we shall consider in this article an equation of the form:
\begin{equation}\label{eq:frac-delay-intro}
dy_t=f(\cz_t^y) dB_t
+ b(\cz_t^y) dt,\quad t\in[0,T],
\end{equation}
where $B$ is a $d$-dimensional fractional Brownian motion with Hurst 
parameter $H>1/2$, 
$f:\cac_1^{\ga}([-h,0];\R^n)\rightarrow \R^{n\times d}$ and
$b:\cac_1^{\ga}([-h,0];\R^n)\rightarrow \R^{n}$ satisfy some suitable
regularity conditions, $\cac_1^\ga$ designates the space of $\ga$-H\"older
continuous functions of one variable (see Section \ref{incr} below) and
$\cz^y_t:[-h,0]\rightarrow\R^n$ is defined by $\cz^y_t(s)=y_{t+s}$. In the 
previous equation,
we also assume that an initial condition
$\xi\in\cac_1^\ga$ is given on the interval $[-h,0]$.
Notice that equation (\ref{eq:frac-delay-intro}) is a slight extension of the
 typical
delay equation which is obtained for some functions $f$ and $b$ of the following
form:
\beq\label{eq:expr-f-weighted-delay}
f:\cac_1^{\ga}([-h,0];\R^n)\rightarrow \R^{n\times d},
\quad\mbox{ with }\quad
f(\cz^y_t)=\si\lp\int_{-h}^0 y_{t+\theta} \, \nu(d\theta) \rp,
\eeq
for a regular enough function $\si$, and a finite measure $\nu$ on $[-h,0]$.
 This special
case of interest will be treated in detail in the sequel. Our considerations
 also include
a function $f$ defined by $f(\cz^y_t)=\si(\cz^y_t(-u_1),\ldots,\cz^y_t(-u_k))$ 
for a given $k\ge 1$,
$0\le u_1<\ldots<u_k\le h$ and a smooth enough function 
$\si:\R^{n\times k}\to\R^{n\times d}$.

\vspace{0.2cm}

The kind of delay stochastic differential system described by 
(\ref{eq:frac-delay-intro}) 
is widely studied when driven by
a standard Brownian motion (see \cite{mohammed} for a nice survey),
 but the results
in the fractional Brownian case are scarce: 
we are only aware of \cite{FR} for the case
$H>1/2$ and $f(\cz^y)=\si(\cz^y(-r)),\ 0\le r\le h$, and the further investigation \cite{FR2} which establishes  a continuity result in terms of the delay $r$. As far as the rough case is concerned, an existence and uniqueness result is given in
 \cite{NNT} for a Hurst parameter $H>1/3$, and \cite{TT} extends this result to $H>1/4$. The current article can be thus
 seen as a
step in the study of processes defined as the solution to fractional delay
 differential
systems, and we shall investigate the behavior of the density of the
 $\R^n$-valued
random variable $y_t$ for a fixed $t\in(0,T]$, where $y$ is the solution to 
(\ref{eq:frac-delay-intro}).
More specifically, we shall prove the following theorem, which can be seen as 
the main
result of the article:
\begin{theorem}\label{thm:smoothness-density-intro}
Consider an equation of the form (\ref{eq:frac-delay-intro}) for an initial condition $\xi$ lying in the space $\cac_1^\ga([-h,0];\R^n)$. Assume $b\equiv 0$, and that $f$ is of the form (\ref{eq:expr-f-weighted-delay}) for a given finite measure $\nu$ on $[-h,0]$ and $\si :\R^n\rightarrow\R^{n\times d}$ a four times differentiable bounded function with bounded derivatives, satisfying the non-degeneracy condition
\begin{equation*}
\si(\eta_1) \si(\eta_2)^* \ge \ep \id_{\R^n},
\quad\mbox{for all}\quad \eta_1,\eta_2\in\R^n.
\end{equation*}
Suppose moreover that $H>H_0$, where $H_0=(7+\sqrt{17})/16\approx 0.6951$.
 Let $t\in(0,T]$ be an arbitrary time, and $y$ be the unique
solution to (\ref{eq:frac-delay-intro}) in $\cac_1^\ka([0,T];\R^n)$, for a given
$1/2<\ka<H$. Then the law of $y_t$ is absolutely continuous with respect to 
Lebesgue
measure in $\R^n$, and its density is a $\cac^\infty$-function.
\end{theorem}
Notice that this kind of result, which has its own interest as a natural step 
in the study of
processes defined by delay systems, is also a useful result when one wants
 to evaluate
the convergence of approximation schemes in the fractional Brownian context.
 We plan to 
report on this possibility in a subsequent communication. The reader may also wonder about our restriction $H>H_0$ above. It will become clear from Remark \ref{rmk:moments-lin} that this assumption is due to the fact that we consider a delay which depends \emph{continuously} on the past. For a discrete type delay of the form $\si(y_t,y_{t-r_1},\ldots,y_{t-r_q})$, with $q\ge 1$ and $r_1<\cdots<r_q\le h$, we shall see at Remark \ref{rmk:discrete-delay} that one can show the smoothness of the density up to $H>1/2$, as for ordinary differential equations. Finally, the case $b\equiv 0$ has been considered here for sake of simplicity, but the extension of our result to a non trivial drift is just a matter of easy additional computations.

\vspace{0.2cm}

Let us say a few words about the strategy we shall follow in order to get our
 Theorem~\ref{thm:smoothness-density-intro}. First of all, as mentioned before, there
 are not too
many results about delay systems governed by a fractional Brownian motion. In particular, equation (\ref{eq:frac-delay-intro}) has never been considered (to the best of our knowledge) with such a general delay dependence. We shall thus first show how to define and solve this differential system, by means of a slight variation of the Young integration theory 
(called algebraic
integration), introduced in \cite{Gu} and also explained in \cite{NNRT}. 
This setting
allows to solve equations like (\ref{eq:frac-delay-intro}) in H\"older spaces 
thanks to 
contraction arguments, in a rather classical way, which will be explained at 
Section \ref{sec:existence-uniqueness}. In fact, observe that our resolution
will be entirely pathwise, and we shall deal with a general 
equation of the form
\begin{equation}\label{eq:general-delay-intro}
dy_t=f(\cz^y_t) dx_t
+ b(\cz^y_t) dt, \quad t\in[0,T],
\end{equation}
for a given  path $x\in \cac_1^{\ga}([0,T];\R^d)$ with $\ga>1/2$, 
where the integral with
respect to $x$ has to be understood in the Young sense \cite{Yo}. 
Furthermore, in equations
like (\ref{eq:general-delay-intro}), the drift term $b(\cz^y)$ is usually 
harmless, but induces
some cumbersome notations. Thus, for sake of simplicity, we shall rather 
deal in the sequel with a reduced delay equation of the type:
$$
y_t=a+\iot f(\cz^y_s)\, dx_s,\quad t\in[0,T].
$$
Once this last equation is properly defined and solved, the 
differentiability of the solution $y_t$
in the Malliavin calculus sense
will be obtained in a pathwise manner, similarly to the case treated in 
\cite{NS}. Finally,
the smoothness Theorem \ref{thm:smoothness-density-intro} will be obtained
mainly by bounding the moments of the Malliavin derivatives of $y$. This will
 be achieved
thanks to a careful analysis and some a priori estimates for equation 
(\ref{eq:frac-delay-intro}).

\smallskip

Here is how our article is structured: Section \ref{sec:alg-intg} is devoted to recall some basic facts about Young integration. We solve, estimate and differentiate a general class of delay equations driven by a Hölder noise at Section \ref{sec:young-delay-eq}. Then at Section \ref{sec:appli-mallia} we apply those general results to fBm and prove our main Theorem \ref{thm:smoothness-density-intro}.

\section{Algebraic Young integration}\label{sec:alg-intg}

The Young integration can be introduced in several ways (convergence of
 Riemann
sums, fractional calculus setting \cite{ZA}).
 We have chosen here to follow the algebraic
approach introduced in \cite{Gu} and developed e.g. in \cite{GT,NNRT}, since
 this 
formalism will help us later in our analysis.

\subsection{Increments}\label{incr}

Let us begin with  the basic algebraic structures which
will allow us to define a pathwise integral with respect to
irregular functions:
first of all,  for an arbitrary real number
$T>0$, a topological vector space $V$ and an integer $k\ge 1$ we denote by
$\cac_k(V)$ (or by $\cac_k([0,T];V)$)
 the set of continuous functions $g : [0,T]^{k} \to V$ such
that $g_{t_1 \cdots t_{k}} = 0$
whenever $t_i = t_{i+1}$ for some $i\le k-1$.
Such a function will be called a
\emph{$(k-1)$-increment}, and we will
set $\cac_*(V)=\cup_{k\ge 1}\cac_k(V)$. An important elementary operator 
is $\der$, which
is defined as follows on $\cac_k(V)$:
\begin{equation}
  \label{eq:coboundary}
\delta : \cac_k(V) \to \cac_{k+1}(V), \qquad
(\delta g)_{t_1 \cdots t_{k+1}} = \sum_{i=1}^{k+1} (-1)^{k-i}
g_{t_1  \cdots \hat t_i \cdots t_{k+1}} ,
\end{equation}
where $\hat t_i$ means that this particular argument is omitted.
A fundamental property of $\der$, which is easily verified,
is that
$\delta \delta = 0$, where $\delta \delta$ is considered as an operator
from $\cac_k(V)$ to $\cac_{k+2}(V)$.
 We will denote $\cz\cac_k(V) = \cac_k(V) \cap \text{Ker}\delta$
and $\cb \cac_k(V) =
\cac_k(V) \cap \text{Im}\delta$.

\vspace{0.2cm}

Some simple examples of actions of $\der$,
which will be the ones we will really use throughout the paper,
 are obtained by letting
$g\in\cac_1(V)$ and $h\in\cac_2(V)$. Then, for any $s,u,t\in\ott$, we have
\begin{equation}
\label{eq:simple_application}
  (\der g)_{st} = g_t - g_s,
\quad\mbox{ and }\quad
(\der h)_{sut} = h_{st}-h_{su}-h_{ut}.
\end{equation}
Furthermore, it is easily checked that
$\cz \cac_{k}(V) = \cb \cac_{k}(V)$ for any $k\ge 1$.
In particular, the following basic property holds:
\begin{lemma}\label{exd}
Let $k\ge 1$ and $h\in \cz\cac_{k+1}(V)$. Then there exists a (non unique)
$f\in\cac_{k}(V)$ such that $h=\der f$.
\end{lemma}

Observe that Lemma \ref{exd} implies that all the elements
$h \in\cac_2(V)$ such that $\der h= 0$ can be written as $h = \der f$
for some (non unique) $f \in \cac_1(V)$. Thus we get a heuristic
interpretation of $\der |_{\cac_2(V)}$:  it measures how much a
given 1-increment  is far from being an  exact increment of a
function, i.e., a finite difference.

\begin{remark}\label{rmk:delta-intg-increment}
Here is a first elementary but important link between these algebraic structures
and integration theory: let $f$ and $g$ be two smooth real valued function on
 $\ott$.
Define then $I\in\cac_2(V)$ by
$$
I_{st}=\ist df_v \int_{s}^{v} dg_w,
\quad\mbox{ for }\quad
s,t\in\ott.
$$
Then, some trivial computations show that
$$
(\der I)_{sut}=[g_u-g_s] [f_t-f_u]=(\der f)_{ut} (\der g)_{su}.
$$
This is a helpful property of the operator $\der$: it transforms iterated 
integrals
into products of increments, and we will be able to take advantage of both
 regularities
of $f$ and $g$ in these products of the form $\der f\, \der g$.
\end{remark}

\vspace{0.2cm}

For sake of simplicity, let us specialize now our setting to the
 case $V=\R^m$ for
an arbitrary $m\ge 1$.
Notice that our future discussions will mainly rely on
$k$-increments with $k \le 2$, for which we will use some
analytical assumptions. Namely,
we measure the size of these increments by H\"older norms
defined in the following way: for $0\le a_1<a_2\le T$ and 
$f \in \cac_2([a_1,a_2];V)$, let
$$
\norm{f}_{\mu ,[a_1,a_2]} =
\sup_{r,t\in[a_1,a_2]}\frac{|f_{rt}|}{|t-r|^\mu},
\quad\mbox{and}\quad
\cac_2^\mu([a_1,a_2];V)=\lcl f \in \cac_2(V);\, 
\norm{f}_{\mu ,[a_1,a_2]}<\infty  \rcl.
$$
Obviously, the usual H\"older spaces 
$\cac_1^\mu([a_1,a_2]; V)$ will be determined
        in the following way: for a continuous function 
$g\in\cac_1([a_1,a_2];V)$,
 we simply set
\begin{equation}\label{def:hnorm-c1}
\|g\|_{\mu,[a_1,a_2]}=\|\der g\|_{\mu ,[a_1,a_2]},
\end{equation}
and we will say that $g\in\cac_1^\mu([a_1,a_2];V)$
 iff $\|g\|_{\mu ,[a_1,a_2]}$ is finite.
Notice that $\|\cdot\|_{\mu ,[a_1,a_2]}$ is only a semi-norm on 
$\cac_1^{\mu}([a_1,a_2];V)$,
but we will generally work on spaces of the type
\begin{equation}\label{def:hold-init}
\cac_{v,a_1,a_2}^\mu(V)=
\lcl g:[a_1,a_2]\to V;\, g_{a_1}=v,\, \|g\|_{\mu ,[a_1,a_2]}<\infty \rcl,
\end{equation}
for a given $v\in V$, or
\begin{equation}\label{eq:def-delayed-norm}
\cac^{\mu}_{\varrho, a_1,a_2}(\R^d):=
\{\zeta\in \cac^{\mu}_{1}
([a_1-h,a_2];\R^d);\zeta=\varrho\ \mbox{\rm on}\ [a_1-h,a_1]\},
\end{equation}
where $0\le a_1<a_2$ and $\varrho\in \cac^{\mu}_1([a_1-h,a_1];\R^d).$
 These last two  spaces are complete metric spaces with the
distance $d_{\mu}(g,f)=\|g-f\|_{\mu}$. More specifically, the metric we shall use
on the space $\cac^{\mu}_{\varrho, a_1,a_2}(\R^d)$ is:
$$
d_{\mu,a_1,a_2}(g,f)\triangleq \|g-f\|_{\mu,[a_1-h,a_2]}.
$$
 In some 
cases we will only write $\cac_k^{\mu}(V)$ instead of 
$\cac_k^{\mu}([a_1,a_2];V)$ when this does not lead to an ambiguity in the
 domain of 
definition of the functions under consideration. 
 For $h \in \cac_3([a_1,a_2];V)$ set in the same way
\begin{eqnarray}
  \label{eq:normOCC2}
  \norm{h}_{\gamma,\rho,[a_1,a_2]} &=& \sup_{s,u,t\in[a_1,a_2]}
\frac{|h_{sut}|}{|u-s|^\gamma |t-u|^\rho}\\
\|h\|_{\mu,[a_1,a_2]} &= &
\inf\left \{\sum_i \|h_i\|_{\rho_i,\mu-\rho_i} ;\, h =
 \sum_i h_i,\, 0 < \rho_i < \mu \right\} ,\nonumber
\end{eqnarray}
where the last infimum is taken over all sequences $\{h_i \in \cac_3(V) \}$
such that $h
= \sum_i h_i$ and for all choices of the numbers $\rho_i \in (0,\mu)$. 
Then  $\|\cdot\|_\mu$ is easily seen to be a norm on $\cac_3([a_1,a_2];V)$,
 and we set
$$
\cac_3^\mu([a_1,a_2];V):=\lcl h\in\cac_3([a_1,a_2];V);\, \|h\|_\mu<\infty \rcl.
$$
Eventually,
let $\cac_3^{1+}([a_1,a_2];V) = \cup_{\mu > 1} \cac_3^\mu([a_1,a_2];V)$,
and remark that the same kind of norms can be considered on the
spaces $\cz \cac_3([a_1,a_2];V)$, leading to the definition of some spaces
$\cz \cac_3^\mu([a_1,a_2];V)$ and $\cz \cac_3^{1+}([a_1,a_2];V)$.

\vspace{0.2cm}

With these notations in mind,
the crucial point in our approach to pathwise integration of irregular
processes is that, under mild smoothness conditions, the operator
$\delta$ can be inverted. This inverse is called $\laa$, and is
defined in the following proposition, whose proof can be found in 
\cite{Gu}. 
\begin{proposition}
\label{prop:Lambda}
Let $0\le a_1<a_2\le T$. Then
there exists a unique linear map $\Lambda: \cz \cac^{1+}_3([a_1,a_2];V)
\to \cac_2^{1+}([a_1,a_2];V)$ such that
$$
\delta \Lambda  = \id_{\cz \cac_3^{1+}([a_1,a_2];V)}.
$$
In other words, for any $h\in\cac^{1+}_3([a_1,a_2];V)$ such that $\der h=0$
there exists a unique $g=\laa(h)\in\cac_2^{1+}([a_1,a_2];V)$ such that 
$\der g=h$.
Furthermore, for any $\mu > 1$,
the map $\laa$ is continuous from $\cz \cac^{\mu}_3([a_1,a_2];V)$
to $\cac_2^{\mu}([a_1,a_2];V)$ and we have
\begin{equation}\label{ineqla}
\|\Lambda h\|_{\mu,[a_1,a_2]} \le \frac{1}{2^\mu-2} \|h\|_{\mu,[a_1,a_2]} ,
\qquad h \in
\cz \cac^{\mu}_3([a_1,a_2];V).
\end{equation}
\end{proposition}

Moreover, the operator $\laa$ can be related to the limit of some Riemann sums,
which gives a second link (after Remark \ref{rmk:delta-intg-increment}) between the previous algebraic developments and some
kind of generalized integration.
\begin{corollary}
\label{cor:integration}
For any 1-increment $g\in\cac_2 (V)$ such that $\der g\in\cac_3^{1+}$,
set
$
\delta f = (\id-\Lambda \delta) g
$.
Then
$$
(\delta f)_{st} = \lim_{|\Pi_{st}| \to 0} \sum_{i=0}^{n-1} g_{t_i\,t_{i+1}},
$$
where the limit is over any partition $\Pi_{st} = \{t_0=s,\dots,
t_n=t\}$ of $[s,t]$, whose mesh tends to zero. Thus, the
1-increment $\delta f$ is the indefinite integral of the 1-increment $g$.
\end{corollary}

\vspace{0.2cm}

\subsection{Young integration}
In this section, we will define a generalized integral $\ist f_u dg_u$ for a 
 $C^{\ka}_1([0,T];\R^{n\times d})$-function $f$, and a 
$\cac_1^{\ga}([0,T];\R^d)$-function
$g$, with $\ka+\ga>1$, by means of the algebraic tools introduced at Section 
\ref{incr}.
To this purpose, we will first assume that $f$ and $g$ are smooth functions,
 in which
case the integral of $f$ with respect to $g$ can be defined in the Riemann 
sense,
and then we will express this integral in terms of the operator $\laa$. This 
will
lead to a natural extension of the notion of integral, which coincides with 
the usual
Young integral. In the sequel, in order to avoid some cumbersome notations,
we will sometimes write $\cj_{st}(f\, dg)$ instead of $\ist f_u dg_u$.

\vspace{0.2cm}

Let us consider then for the moment two smooth functions $f$ and $g$
defined on $\ott$. One can
write, thanks to some elementary algebraic manipulations, that:
\begin{equation}\label{eq:def-intg-smooth-case}
\cj_{st}(f\, dg)\equiv \ist f_u \, dg_u
= f_s (\der g)_{st} + \ist (\der f)_{su} \, dg_u
=f_s (\der g)_{st} + \cj_{st}(\der f\, dg).
\end{equation}
Let us analyze now the term $\cj(\der f\, dg)$, which is an element of 
$\cac_2(\R^n)$. 
Invoking Remark \ref{rmk:delta-intg-increment}, it is easily seen that, for 
$s,u,t\in\ott$,
$$
h_{sut}\equiv \lc \der\lp \cj(\der f\, dg) \rp \rc_{sut} = (\der f)_{su} 
(\der g)_{ut}.
$$
The increment $h$ is thus an element of $\cac_3(\R^n)$ satisfying $\der h=0$ (recall that $\der\der=0$).
 Let us estimate
now the regularity of $h$: if $f\in C^{\ka}_1([0,T];\R^{n\times d})$ and $g\in
 \cac_1^{\ga}([0,T];\R^d)$,
 from the definition 
(\ref{eq:normOCC2}), it is readily checked that $h\in\cac_3^{\ga+\ka}(\R^n)$. 
Hence
$h\in\cz\cac_3^{\ga+\ka}(\R^n)$, and if $\ka+\ga>1$ (which is the case if $f$ 
and $g$ are regular),
Proposition \ref{prop:Lambda} yields that $\cj(\der f\, dg)$
can also be expressed as
$$
\cj(\der f\, dg) = \laa(h) = \laa\lp  \der f \, \der g\rp,
$$
and thus, plugging this identity into (\ref{eq:def-intg-smooth-case}), we get:
\begin{equation}\label{eq:def2-intg-smooth-case}
\cj_{st}(f\, dg)
=f_s (\der g)_{st} + \laa_{st}\lp  \der f \, \der g\rp.
\end{equation}
Now we can see that the right hand side of the last equality is rigorously 
defined 
whenever $f\in C^{\ka}_1([0,T];\R^{n\times d})$, $g\in \cac_1^{\ga}([0,T];\R^d)$,
 and this is the definition we will use in order to 
extend the notion of integral:
\begin{theorem}\label{prop-cj}
Let $f \in \cac^{\ka}_1([0,T];\R^{n\times d})$ and $g
\in \cac_1^{\ga}([0,T];\R^d)$, with $\ka+\ga>1$. Set
\begin{equation}\label{eq:def-young-intg}
\cj_{st}( f\, dg)= f_s (\der g)_{st} + \laa_{st}\lp  \der f \, \der g\rp.
\end{equation}
Then
\begin{itemize}
\item[(1)]
Whenever $f$ and $g$ are smooth function, $\cj_{st}( f\, dg)$ coincides with the
usual Riemann integral.
\item[(2)]
The generalized integral $\cj( f\, dg)$ satisfies:
$$
\lln \cj_{st}( f\, dg) \rrn \le  \|f\|_{\infty} \|g\|_{\ga} |t-s|^\ga
+ c_{\ga,\ka} \|f\|_{\ka} \|g\|_{\ga} |t-s|^{\ga+\ka},
$$
for a constant $c_{\ga,\ka}$ whose exact value is $(2^{\ga+\ka}-1)^{-1}$.
\item[(3)]
We have
$$
\cj_{st}( f\, dg) = \lim_{|\Pi_{st}| \to 0} \sum_{i=0}^{n-1}  f_{t_i} \, 
\der g_{t_i\,t_{i+1}},
$$
where the limit is over any partition $\Pi_{st} = \{t_0=s,\dots,
t_n=t\}$ of $[s,t]$, whose mesh tends to zero. In particular, 
$\cj_{st}( f\, dg)$ coincides
with the Young integral as defined in \cite{Yo}.
\end{itemize}
\end{theorem}

\begin{proof}
The first claim is just what we proved at equation 
(\ref{eq:def2-intg-smooth-case}).
The second assertion follows directly from the definition 
(\ref{eq:def-young-intg}) and 
the inequality (\ref{ineqla}) concerning the operator $\laa$. Finally, 
our third property
is a direct consequence of Corollary \ref{cor:integration} and the fact that
$\der(f\, \der g)=-\der f \der g$, which means that
$$
\cj( f\, dg) = \lc \id-\laa\der \rc(f\, \der g).
$$

\end{proof}

A Fubini type theorem for Young's integral will be needed
 in the last section of this paper. Its proof below is a good example of the 
importance of Proposition \ref{prop:Lambda} and Theorem \ref{prop-cj}.

\begin{proposition}\label{prop:25f}
Assume that $\ga>\la>1/2$. Let $f$ and $g$ be two functions in
$\cac_1^\ga([0,T]:\R)$ and $h:\{(t,s)\in [0,T]^2; 0\le s\le t\le T\}
\rightarrow\R$ a function such that $h(\cdot,t)$ (resp. $h(t,\cdot)$)
belongs to $\cac^\la_1([t,T];\R)$ (resp. $\cac^\la_1([0,t];\R)$) uniformly
in $t\in[0,T]$, and
\begin{equation}\label{eq:uhc}
\| h(r_1,\cdot)-h(r_2,\cdot)\|_{\la,[0,r_1\wedge r_2]}\le
C|r_1-r_2|^\la.\end{equation}
Then 
\begin{equation}\label{eq:ft}
\int_s^t\int_s^r h(r,u)dg_udf_r=
\int_s^t\int_u^t h(r,u)df_r dg_u,\quad 0\le s\le t\le T.\end{equation}
\end{proposition}

\begin{proof}
Fix $s,t\in[0,T]$, with $s<t$, and divide the proof in several steps.

\smallskip

\noindent\textit{Step 1.} Here we see that $\int_s^t\int_s^r h(r,u)
dg_udf_r$ is well-defined. Note that we only need to show that
$\int_s^\cdot h(\cdot,u)dg_u$ belongs to  $\cac^\la_1([s,T];\R)$ due
to Theorem \ref{prop-cj}.

Let $r_1,r_2\in[s,t],\ r_1<r_2$, then Theorem \ref{prop-cj}.(2) gives
\begin{eqnarray*}
\lefteqn{\left|\int_s^{r_2} h(r_2,u)dg_u-\int_s^{r_1} h(r_1,u)dg_u\right|}\\
&\le& \left|\int_s^{r_1}(  h(r_2,u)-h(r_1,u))dg_u\right|
+\left|\int_{r_1}^{r_2} h(r_2,u)dg_u\right|\\
&\le&\|g\|_\ga\left( \|h(r_2,\cdot)-h(r_1,\cdot)\|_{\infty,[0,r_1]}
(r_1-s)^\ga +c_{\ga,\la} \|h(r_2,\cdot)-h(r_1,\cdot)\|_{\la,[0,r_1]}
(r_1-s)^{\ga+\la}\right)\\
&&+\|g\|_\ga\left( \|h(r_2,\cdot)\|_{\infty,[0,r_2]}
(r_2-r_1)^\ga +c_{\ga,\la} \|h(r_2,\cdot)\|_{\la,[0,r_2]}
(r_2-r_1)^{\ga+\la}\right).
\end{eqnarray*}
Hence (\ref{eq:uhc}) implies our claim. The definition of $\int_s^t\int_u^t h(r,u)df_r dg_u$ follows along the same lines.

\smallskip

\noindent\textit{Step 2.}   Let $\Pi_{st}=\{t_0=s,\ldots, t_n=t\}$ be a partition
of the interval $[s,t]$. Then, according to Proposition \ref{prop-cj}, for any $v\in[0,t)$ we have
\begin{equation}\label{eq:fubi1}
\int_s^v h(t,u)dg_u=\lim_{|\Pi_{st}|\rightarrow 0}\sum_{i=0}^{n-1}
h(t,t_i) \, (\der g)_{t_i\wedge v, t_{i+1}\wedge v}.
\end{equation}
Our assumption (\ref{eq:uhc}) allows now to take limits in the equation above, so that we obtain, for any $0\le s<t\le T$,
\begin{equation}\label{eq:fu4}
q_{st}^{1}:=\int_s^t h(t,u)dg_u
=\lim_{|\Pi_{st}|\rightarrow 0}\sum_{i=0}^{n-1} h(t,t_i) \, \der g_{t_i,t_{i+1}}:=q_{st}^{2}.
\end{equation}
In order to see that the relation above holds in $\cac_2^{\la}(\ott;\R)$, it is now enough to check that both $q^1$ and $q^2$ in  (\ref{eq:fu4}) are elements of $\cac_2^{\la}(\ott;\R)$.

\smallskip

However, the fact that $q^1\in\cac_2^{\la}(\ott;\R)$ can be proved along the same lines as in Step~1. The assertion $q^2\in\cac_2^{\la}(\ott;\R)$ can be proved by observing that the limit defining $q_{st}^{2}$ do not depend on the sequence of partitions under consideration. In particular, consider the sequence $(\pi^n)_n$ of dyadic partitions of $[0,T]$, that is
$$
\pi^n=\{ 0=t_0^n \leq t_1^n \leq \dots \leq t_{2^n}^n =T\}, \ \mbox{with} \ t_i^n=\frac{i\, T}{2^n},
$$
and set, for all $s,t\in [0,T]$, $\pi_{st}^{n}=\pi^{n}\cap (s,t)$. Then
$
q_{st}^{2}=
\lim_{n\to\infty}\sum_{t_i\in\pi_{st}^{n}} h(t,t_i^{n}) \, \der g_{t_i^{n},t_{i+1}^{n}}
$
for all $0\le s<t\le T$,
and the same kind of arguments as in \cite[Theorem 2.2]{DT2} yield our claim $q^2\in\cac_2^{\la}(\ott;\R)$. We have thus proved that (\ref{eq:fu4}) holds in $\cac_2^{\la}(\ott;\R)$.

\smallskip

\noindent\textit{Step 3.} From Proposition \ref{prop:Lambda}, Step 2 and 
(\ref{eq:def-young-intg}) we have 
\begin{align*}
&\int_s^t\int_s^r h(r,u)dg_udf_r=\lim_{|\Pi_{st}|\rightarrow 0}
\int_s^t\left(\sum_{i=0}^{n-1} h(r,t_i)(g_{t_{i+1}\wedge r}-g_{t_{i}\wedge r})\right)
df_r \\
&
=\lim_{|\Pi_{st}|\rightarrow 0}\sum_{i=0}^{n-1} \int_{t_i}^t
h(r,t_{i}) \left(g_{t_{i+1}\wedge r}-g_{t_i}\right)\, df_r \\
&=\lim_{|\Pi_{st}|\rightarrow 0}\sum_{i=0}^{n-1}\lc \left(\int_{t_i}^t
h(r,t_i)df_r\right)\left(g_{t_{i+1}}-g_{t_i}\right)
+ \int_{t_i}^{t_{i+1}} h(r,t_{i}) \left(g_{t_{i+1}\wedge r}-g_{t_{i+1}}\right)\, df_r \rc \\
\end{align*}
Moreover, thanks to the Hölder properties of $f$ and $g$, we have
$$\sum_{i=0}^{n-1}\left|\int_{t_i}^{t_{i+1}}h(r,t_i)(g_{r}-g_{t_{i}})df_r\right|
\le C\sum_{i=0}^{n-1}(t_{i+1}-t_{i})^{\ga+\la}\rightarrow 0$$
as $|\Pi_{st}|\rightarrow 0$, and thus
$$
\int_s^t\int_s^r h(r,u)dg_udf_r=
\lim_{|\Pi_{st}|\rightarrow 0}\sum_{i=0}^{n-1} \left(\int_{t_i}^t
h(r,t_i)df_r\right)\left(g_{t_{i+1}}-g_{t_i}\right).
$$
Consequently, Step 2 and Theorem \ref{prop-cj}
imply that (\ref{eq:ft}) is satisfied and therefore the proof is complete.
\end{proof}

The following integration by parts  and It\^o's formulas 
will be also needed in the last
part of this paper.

\begin{proposition}\label{prop:ipf} Let $f$ and $g$ be two functions in
$\cac_1^\ga([0,T];\R)$, with $\ga>1/2$. Then 
$$f_tg_t=f_0g_0+\int_0^tf_udg_u+\int_0^tg_udf_u,\quad t\in[0,T].$$
\end{proposition}
\begin{proof}Set $q_t:=f_tg_t-\int_0^tf_udg_u-\int_0^tg_udf_u$, $t\in[0,T]$.
It is easy to see that 
this  funcion belongs to $\cac_1^{2\ga}([0,T];\R)$ because of the equalities
$$f_tg_t-f_sg_s=f_s(\delta g)_{st}+g_s(\delta f)_{st}+(\delta g)_{st}
(\delta f)_{st}$$
and
$$\int_s^tf_udg_u+\int_s^t g_udf_u=f_s(\delta g)_{st}+g_s(\delta f)_{st}
+\laa_{st}(\delta f\delta g)+\laa_{st}(\delta g\delta f),$$
which follows from 
(\ref{eq:def-young-intg}). Now, since $q\in\cac_1^{2\ga}([0,T];\R)$, with $2\ga>1$, $q$ is a constant function. Otherwise stated,
$q_t=q_0=f_0g_0$. Therefore the announced result is true.
\end{proof}
\begin{proposition} \label{prop:ifh}
Let $g$ and $h$ be in $\cac_1^\ga([0,T],\R)$ and
$f\in\cac^2_b(\R)$. Also let $x_t=x_0+\int_0^t g_sdh_s$, $t\in[0,T]$. Then
$$f(x_t)=f(x_0)+\int_0^tf'(x_u)g_udh_u,\quad t\in[0,T].$$
\end{proposition}
\begin{proof} Proceeding as in the proof of Proposition \ref{prop:ipf}
and using the mean value theorem, we can show that 
$$q_t=f(x_t)-\int_0^t f'(x_s)g_sdh_s,\quad t\in[0,T],$$
is a $2\ga$-H\"older continuous function. Therefore the result holds.
\end{proof}

\begin{remark} 
Proposition \ref{prop:ifh} has been proven in \cite{ZA}
using Riemann sums.
\end{remark}

\section{Young delay equation}\label{sec:young-delay-eq}
Recall first that we wish to consider a differential equation of the form:
\begin{eqnarray} \label{deleq}
 y_t&=&\xi_0+\int_0^t f(\cz^y_u)\,dx_u, \quad t\in[0,T],\\
\cz^y_0&=&\xi .\nonumber
\end{eqnarray}
In the previous equation, the integral  has to be interpreted in the Young 
sense of 
(\ref{eq:def-young-intg}), the initial condition $\xi$ is an element of
$\cac^{\ga}_1([-h,0];\R^n)$,  the driving noise $x$ is in
$\cac^{\ga}_1([0,T];\R^d)$, with $\ga>1/2$. We seek a solution $y$ in the
space $\cac^{\la}_{\xi,0,T}(\R^n)$ for $1/2<\la<\ga$, and $f$ is a given function
 $f:\cac_1^{\la}([-h,0];\R^n)\rightarrow \R^{n\times d}$. In this section, we 
shall solve
 equation (\ref{deleq}) thanks to a contraction argument, and then study its
 differentiability with respect to the driving noise $x$. Of course, the main 
application
 we have in mind is the case where $x$ is a $d$-dimensional fractional
 Brownian motion,
 and this particular case will be considered at Section \ref{sec:appli-mallia}.
 
\vspace{0.2cm}

\subsection{Existence and uniqueness of the solution}\label{sec:existence-uniqueness}
In order to solve equation (\ref{deleq}), some smoothness and boundedness 
assumptions
have to be made on our coefficient $f$. In fact, we shall rely on the 
following hypothesis:
\begin{hypothesis}\label{hyp:H1}
There exist a positive constant $M$ 
and $\la\in(1/2,\ga)$ such that
$$
|f(\zeta)|\le M,
\quad\mbox{ and }\quad
|f(\zeta_2)-f(\zeta_1)|\le M \sup_{\theta\in[-h,0]}|\zeta_2(\theta)-\zeta_1
(\theta)|
$$
uniformly in $\zeta,\zeta_1,\zeta_2\in \cac^{\la}_1([-h,0];\R^n)$.
\end{hypothesis}
\noindent
Actually we will assume that $f$ satisfies a stronger
Lipschitz type  hypothesis on the space
$\cac^{\la}_1(\R^n)$. Let us state first a preliminary result before we come to
 this second
assumption:
\begin{lemma}\label{integrand}
Let $a=(a_1,a_2)$, with $0\le a_1<a_2\le T$, 
let also $Z\in\cac^{\la}_1([a_1-h,a_2];\R^n)$ and set
$$\left[\UU^{(a)} Z\right]_s=f(\cz^Z_s),\quad s\in[a_1,a_2].$$
Then Hypothesis \ref{hyp:H1} implies that  $\UU^{(a)}$ is a
map from  
$\cac^{\la}_{1}([a_1-h,a_2];\R^n)$ into
$
\cac^{\la}_1([a_1,a_2];$ $\R^{n\times d})$, satisfying:
$$
\left\| \UU^{(a)} Z \right\|_{\la,[a_1,a_2]} \le  M \,
 \left\|  Z \right\|_{\la,[a_1-h,a_2]}.
$$
 \end{lemma}

\begin{proof} 
 The proof of this result is an immediate consequence of the definition
(\ref{def:hnorm-c1}) of H\"older's norms on $\cac_1$ and Hypothesis
\ref{hyp:H1}.

\end{proof}

With this preliminary result in hand, we can now introduce our second 
hypothesis on the coefficient 
$f$.
\begin{hypothesis}\label{hyp:H2}
Taking up the notations of Hypothesis \ref{hyp:H1}, 
consider an initial condition
$\rho\in\cac_{1}^{\la}([a_1-h,a_1])$. We assume that,
for any $N\ge 1$, there is a positive constant $c_N$ such that:
$$\|\UU^{(a)}(Z_1)-\UU^{(a)}(Z_2)\|_{\la,[a_1,a_2]}
\le c_N \|Z_1-Z_2\|_{\la,[a_1-h,a_2]},$$
for all $0\le a_1\le a_2\le T$ and $Z_1,Z_2\in
\cac^{\la}_{\rho,a_1,a_2}(\R^n)$, satisfying 
$$
\max\lcl \|Z_1\|_{\la,[a_1-h,a_2]};\,  \|Z_2\|_{\la,[a_1-h,a_2]} \rcl \le N,
$$
where $\la$ is given in Hypothesis \ref{hyp:H1}.
\end{hypothesis}
\noindent
Observe that Hypothesis \ref{hyp:H2} holds in particular if, for
$\la>0$, the map $\UU^{(a)}$ admits a derivative which is locally bounded,
 uniformly in $a\in[0,T]$. 

\vspace{0.2cm}

Now that we have stated our main assumptions, the following theorem  is the
 main result of this section.
\begin{theorem}\label{exiuniq} Under 
Hypotheses \ref{hyp:H1} and \ref{hyp:H2}, 
 the delay equation (\ref{deleq}) has a 
unique solution in $\cac^{\la}_{\xi, 0,T}(\R^n)$.
\end{theorem}

Before giving the proof of this theorem, we establish and auxiliary result.
 This will
be helpful in order to get the existence of an invariant ball under the 
contracting map which 
gives raise to the solution of our equation.

\begin{lemma}\label{ineq-int}
 Let $x\in\cac_1^\ga([a_1,a_2];\R^d)$
 with $\ga>1/2$ and $0\le a_1<a_2$, $\la\in(1/2,\ga)$ and
$v\in\R^n$. Set $a=(a_1,a_2)$, recall notation (\ref{def:hold-init}) and define
  $\VV^{(a)}:\cac_1^{\la}([a_1,a_2];\R^{n\times d})
\rightarrow \cac_{v,a_1,a_2}^{\la}(\R^n)$ by:
$$\left[\VV^{(a)} Z\right]_s= v+\cj_{a_1s}(Z\, dx),\quad s\in[a_1,a_2],$$
where $\cj_{a_1s}(Z\, dx)$ stands for the Young integral defined by 
(\ref{eq:def-young-intg}). 
Then 
$$\|\VV^{(a)} Z\|_{\la,[a_1,a_2]}\le \|x\|_{\ga}
\left(\|Z\|_{\infty,[a_1,a_2]}(a_2-a_1)^{\ga-\la}+
c_{\la+\ga}\| Z\|_{\la,[a_1,a_2]}(a_2-a_1)^{\ga}\right),$$
with $c_{\la+\ga}=(2^{\la+\ga}-2)^{-1}$.
\end{lemma}

\begin{proof}
 Let $a_1\le s\le t\le T$. Then Theorem \ref{prop-cj} point (3) implies that
$$\left[\VV^{(a)} Z\right]_t- \left[\VV^{(a)}
 Z\right]_s=\cj_{st}(Z\, dx).$$
 Our claim is then a direct consequence of Theorem \ref{prop-cj} point (2)
 and of the definition  (\ref{def:hnorm-c1}).
 
\end{proof}

\begin{proof}[Proof of Theorem \ref{exiuniq}:]
This proof is divided in several steps.

\vspace{0.2cm}

\noindent
{\it Step 1: Existence of invariant balls.}
Let us first consider an interval of the form $[0,\ep]$, which means that, 
when we include
the delay of the equation, we shall consider processes defined on $[-h,\ep]$.
 More specifically,
let us recall that the spaces $\cac^{\la}_{\xi,0,\ep}(\R^n)$
 have been defined by relation 
(\ref{eq:def-delayed-norm}). Then we consider a map 
$\gaa:\cac^{\la}_{\xi,0,\ep}\to\cac^{\la}_{\xi,0,\ep}$, where
we have set $\cac_{\xi,0,\ep}^{\la}=\cac_{\xi,0,\ep}^{\la}(\R^n)$ for notational 
sake,
defined in the following way: if $z\in\cac_{\xi,0,\ep}^{\la}$, then $\gaa(z)=\hz$,
 where
$\hz_t=\xi_t$ for $t\in[-h,0]$, and:
\begin{equation}\label{eq:def-der-hat-z}
(\der\hz)_{st}=\cj_{st}(Z\, dx),
\quad\mbox{ with }\quad Z_u=f(\cz^z_u),
\quad\mbox{ for }\quad s,t\in[0,\ep].
\end{equation}
We shall now look for an invariant ball in the space
 $\cac^{\la}_{\xi,0,\ep}$ for
 the map $\gaa$. 

\vspace{0.2cm}

So let us pick an element $z$, such that $\|z\|_{\la,[-h,\ep]}\le N_1$ and
set $\gaa(z)=\hz$. On $[-h,0]$, we have $\hz=\xi$, and hence 
$\|\der\hz\|_{\la,[-h,0]}=\|\der\xi\|_{\la,[-h,0]}\equiv N_\xi$. We shall thus 
choose
$N_1\ge 2N_\xi$.

\vspace{0.2cm}

On $[0,\ep]$, we have now, invoking Lemma \ref{ineq-int}:
\beq\label{eq:bnd-der-hz}
\|\der\hz\|_{\la,[0,\ep]}
\le
\|Z\|_{\infty} \|x\|_\ga \ep^{\ga-\la} + c_{\ga,\la} \|Z\|_{\la,[0,\ep]} 
\|x\|_\ga \ep^{\ga}.
\eeq
Furthermore, according to Hypothesis \ref{hyp:H1}, we have $\|Z\|_{\infty}\le M$
and thanks to Lemma \ref{integrand}, we also have 
$\|Z\|_{\la,[0,\ep]}\le M \,\|z\|_{\la,[-h,\ep]}\le M \,N_1$, by assumption. Then 
we can recast
the previous inequality into:
\begin{equation}\label{eq:ineq1-der-hz}
\|\der\hz\|_{\la,[0,\ep]}
\le
M \, \|x\|_\ga \ep^{\ga-\la} \lc 1+  c_{\ga,\la} N_1 \ep^{\la}\rc .
\end{equation}
Let us choose now $\ep$ and $N_1$ in the following manner 
(notice that $\ep$ does
{\bf not} depend on the initial condition $\xi$):
\begin{equation}\label{eq:cdt-ep-N_1}
\ep=\lc  4 M c_{\ga,\la} \|x\|_\ga \rc^{-1/\ga}\wedge 1,
\quad\mbox{ and }\quad
N_1\ge 4 M  \|x\|_\ga.
\end{equation}
With this choice of $\ep,N_1$, inequality (\ref{eq:ineq1-der-hz}) becomes
$\|\der\hz\|_{\la,[0,\ep]}\le N_1/2$. Summarizing the considerations above, 
we have thus found that:
\begin{multline}\label{eq:inv-ball1}
\ep=\lc  4 M c_{\ga,\la} \|x\|_\ga \rc^{-1/\ga}\wedge 1, \,
N_1\ge\sup\lcl 2N_\xi;\, 4 M  \|x\|_\ga \rcl \\
\Longrightarrow
\sup\lcl \|\der\hz\|_{\la,[-h,0]};\, \|\der\hz\|_{\la,[0,\ep]}\rcl\le 
\frac{N_1}{2}.
\end{multline}

\vspace{0.2cm}

Consider now $s<t$, with $s\in[-h,0]$ and $t\in[0,\ep]$. Then, owing to the 
previous
relation, we have:
$$
|(\der \hz)_{st}|\le |(\der \hz)_{s0}| + |(\der \hz)_{0t}| \le \frac{N_1}{2} 
\lp s^\la+t^\la \rp
\le N_1 |t-s|^\la,
$$
which, together with the last inequality, proves that $B(0,N_1)$ in 
$\cac_{\xi,0,\ep}^{\la}$
is left invariant by $\gaa$, under the assumptions of (\ref{eq:inv-ball1}). 

\vspace{0.2cm}

Assume now that we have been able to produce a solution $y^{(1)}$ to 
equation (\ref{deleq}) on the interval $[-h,\ep]$. We try now to iterate the
 invariant
ball argument on $[\ep-h;2\ep]$. The arguments above go through with very little
changes: we are now working on delayed
H\"older spaces of the form $\cac_{y^{(1)},\ep,2\ep}^{\la}$, 
and the map $\gaa$ is defined by $\gaa(z)=\hz$, with $\hz=y^{(1)}$ on
$[\ep-h;\ep]$, and $\der\hz$ having the same expression as in 
(\ref{eq:def-der-hat-z})
on $[\ep,2\ep]$. We wish to find a ball $B(0,N_2)$ in  
$\cac_{y^{(1)},\ep,2\ep}^{\la}$,
left invariant by the map $\gaa$.
With the same computations as for the interval $[-h,\ep]$, the
assumptions of inequality (\ref{eq:inv-ball1}) become:
$$
\ep=\lc  4 M c_{\ga,\la} \|x\|_\ga \rc^{-1/\ga}\wedge 1, \,
N_2\ge\sup\lcl 2N_{y^{(1)}};\, 4 M  \|x\|_\ga \rcl.
$$
Notice again that we are able to choose here the same $\ep$ as before, by 
changing
$N_1$ into $N_2$ according to the value of $\|y^{(1)}\|_{\la,[\ep-h,\ep]}$. It is
 now
readily checked that $B(0,N_2)$ is invariant under $\gaa$, and this calculation
is also easily repeated on any interval $[k\ep-h,(k+1)\ep]$ for any $k\ge 0$,
 until
the whole interval $\ott$ is covered.

\vspace{0.2cm}

\noindent
{\it Step 2: Fixed point argument.}
We shall suppose here that we have been able to construct the unique solution 
$y$
to (\ref{deleq}) on $[-h;l\ep]$, and we shall build the fixed point argument on 
$[l\ep-h;(l+1)\ep]$. On the latter interval, the initial condition of the 
paths we shall
consider is $\xi^{l,1}\equiv y$ on
$[l\ep-h;l\ep]$. If $\gaa$ is the map defined on
$\cac_{\xi^{l,1},l\ep,(l+1)\ep}^{\la}$ by (\ref{eq:def-der-hat-z}), then we know
 that
$B(0,N_{l+1})$ is invariant by $\gaa$.

\vspace{0.2cm}

In order to settle our fixed point argument, we shall first consider an
 interval of the 
form $[l\ep-h;l\ep+\eta]$, for a parameter $0<\eta\le\ep$ to be determined. On
$\cac_{\xi^{l,1},l\ep,l\ep+\eta}^{\la}$, we define a map, called again $\gaa$,
according to  (\ref{eq:def-der-hat-z}). Pick then two functions 
$z^1,z^2\in\cac_{\xi^{l,1},l\ep,l\ep+\eta}^{\la}$, set $\hz^i=\gaa(z^i)$ for 
$i=1,2$
and $\zeta=\hz^2-\hz^1$. Then $\zeta\in\cac_{0,l\ep,l\ep+\eta}^{\la}$, and if
$l\ep\le s< t\le l\ep+\eta$, we have
$$
(\der\zeta)_{st}=\cj_{st}\lp  (Z^{2}-Z^{1}) \, dx \rp,
\quad\mbox{ where }\quad
Z^i=f(\cz^{z^i}).
$$
Thus, just like in (\ref{eq:bnd-der-hz}), we have:
$$
\|\der\zeta\|_{\la,[l\ep-h,l\ep+\eta]} \le
\|Z^1-Z^2\|_{\infty,[l\ep,l\ep+\eta]} \|x\|_{\ga} \eta^{\ga-\la}
+c_{\ga,\la} \|Z^1-Z^2\|_{\la,[l\ep,l\ep+\eta]} \|x\|_{\ga} \eta^{\ga}.
$$
Furthermore, 
$\|Z^1-Z^2\|_{\infty,[l\ep,l\ep+\eta]}
\le\|Z^1-Z^2\|_{\la,[l\ep,l\ep+\eta]}\,\eta^{\la}$. Hence,
$$
\|\der\zeta\|_{\la,[l\ep-h,l\ep+\eta]} \le
(1+c_{\ga,\la})\, \|Z^1-Z^2\|_{\la,[l\ep,l\ep+\eta]} \,\|x\|_{\ga} \,\eta^{\ga}.
$$
We also have $Z^1-Z^2=f(\cz^{z^1})-f(\cz^{z^2})$, and thanks to Hypothesis
\ref{hyp:H2}, we obtain:
$$
\|\der\zeta\|_{\la,[l\ep-h,l\ep+\eta]} \le
(1+c_{\ga,\la})\, \|x\|_{\ga} \, c_{N_{l+1}}\,\eta^{\ga}\, 
\|z^1-z^2\|_{\la,[l\ep-h,l\ep+\eta]}.
$$
Therefore,
 we are able to apply the fixed point argument in the usual way as soon as
$$
(1+c_{\ga,\la})\, c_{N_{l+1}}\, \|x\|_{\ga} \,\eta^{\ga}  \le \frac12,
\quad\mbox{ or }\quad
\eta= \lc 2  (1+c_{\ga,\la})\, c_{N_{l+1}}\, \|x\|_{\ga} \rc^{-1/\ga}\wedge \ep.
$$
With this value of $\eta$, we are thus able to get a unique solution to
(\ref{deleq}) on $[l\ep-h;l\ep+\eta]$.

\vspace{0.2cm}

Let us proceed now to the case of $[l\ep+\eta-h,l\ep+2\eta]$. The arguments are
roughly the same as in the previous case, but one has to be careful about the
change in the initial condition. In fact, the initial condition here should be
$\xi^{l,2}\equiv y$ on $[l\ep+\eta-h,l\ep+\eta]$.
 However, we can also choose to extend
this initial condition backward, and set it as $\xi^{l,2}\equiv y$ on
$[l\ep-h,l\ep+\eta]$.
We then define the usual map $\gaa$ as in (\ref{eq:def-der-hat-z}), and we have
to prove that $B(0,N_{l+1})$ is left invariant by $\gaa$. To this purpose,
take $z\in\cac_{\xi^{l,2},l\ep+\eta,l\ep+2\eta}^{\la}$ in $B(0,N_{l+1})$, and set
 $\hz=\gaa(z)$.
Observe then that, for any $t\in[l\ep+\eta,l\ep+2\eta]$, we have
$$
\hz_t=\xi_{l\ep+\eta}^{2}+\int_{l\ep+\eta}^{t} f(\cz^z_u) \, dx_u
= \xi_{l\ep}^{1}+\int_{l\ep}^{l\ep+\eta} f(\cz^y_u) \, dx_u+
\int_{l\ep+\eta}^{t} f(\cz^z_u) \,
 dx_u
= \xi_{l\ep}^{1}+\int_{l\ep}^{t} f(\cz^z_u) \, dx_u,
$$
where we have used the fact that $\xi^{l,2}\equiv y$ on $[l\ep-h,l\ep+\eta]$
 solves 
(\ref{deleq}).
It is now easily seen that $\hz$ is in $B(0,N_{l+1})$, and this allows to 
settle our
fixed point argument as in the previous case, with the same interval
length $\eta$. This step can now be iterated until the whole interval
 $[l\ep;(l+1)\ep]$
is covered.

\end{proof}

\subsection{Moments of the solution}

The moments of the solution to (\ref{deleq}) can be bounded in the following
 way:
\begin{proposition}\label{prop:moments-solution-delay-eq}
Under the same assumptions as in Theorem \ref{exiuniq}, let $y$ be the solution of equation (\ref{deleq}) on the interval $\ott$, 
with an initial condition
$\xi\in\cac_1^\la([-h,0];\R^n)$. Then there exists a strictly positive constant 
$c=c(\ga,\la,M,T)$ such that
$$
\|y\|_{\la,[-h,T]} \le c \max\lc \|\xi\|_\la , \|x\|_\ga^{\la/(\ga+\la-1)},
\|x\|_\ga\rc.
$$
\end{proposition}

\begin{proof}
From the proof of Theorem \ref{exiuniq}, we know that $\|y\|_{\la,[-h,T]}$ is
 finite.
Let us assume that this quantity is equal to $K$, 
and let us find an estimate on
$K$. One can begin with a small interval, which will be called again $[0,\ep]$,
though it won't
be the same interval as in the proof of Theorem \ref{exiuniq}. In any case,
taking into account that $y$ solves equation (\ref{deleq}),
we obtain similarly to (\ref{eq:bnd-der-hz}):
\begin{eqnarray}\label{numagre}
\|\der y\|_{\la,[0,\ep]}&\le&
M \, \|x\|_\ga \ep^{\ga-\la} + c_{\ga,\la} \, M\, \|\der y\|_{\la,[-h,\ep]} 
\, \|x\|_\ga \ep^{\ga} \nonumber \\
&\le& M \, \|x\|_\ga \ep^{\ga-\la} + c_{\ga,\la} \, M\, K \, 
\|x\|_\ga \ep^{\ga} \equiv g(\ep,K).
\end{eqnarray}
Along the same line, for any $k\le[T/\ep]$, we have
$$
\|\der y\|_{\la,[k\ep,(k+1)\ep]} \le g(\ep,K).
$$
Take now $s,t\in\ott$ such that $i\ep\le s<(i+1)\ep\le j\ep\le t <(j+1)\ep$.
 Set also
$t_i=s$, $t_k=k\ep$ for $i+1\le k\le j$, and $t_{j+1}=t$. 
Then
$$
\lln (\der y)_{st}\rrn
= \lln \sum_{k=i}^{j}(\der y)_{t_k t_{k+1}}\rrn
\le g(\ep,K) \sum_{k=i}^{j} (t_{k+1}-t_{k})^{\la}
\le g(\ep,K) (j-i+1)^{1-\la} (t-s)^\la,
$$
where we have used the fact that $r\mapsto r^\la$ is a concave function.
Note that the indices $i,j$ above satisfy
$(j-i+1)\le 2T/\ep$. Plugging this into the last series of inequalities, 
we end up with
$$
\|\der y\|_{\la,[0,T]} \le \frac{g(\ep,K)(2T)^{1-\la}}{\ep^{1-\la}}
= \lc \frac{M \, \|x\|_\ga}{\ep^{1-\ga}} + c_{\ga,\la} \, M\, K \, \|x\|_\ga 
\ep^{\ga+\la-1} \rc (2T)^{1-\la}.
$$
Thus the parameters $K$ and $\ep$  satisfy the relation:
\beq\label{eq:rel-epsilon-K}
K\le 
\lc \frac{M \, \|x\|_\ga}{\ep^{1-\ga}} + c_{\ga,\la} \, M\, K \, 
\|x\|_\ga \ep^{\ga+\la-1} \rc (2T)^{1-\la}+\|\xi\|_{\la},
\eeq
In order to solve (\ref{eq:rel-epsilon-K}), choose $\ep$ such that
$$
c_{\ga,\la} \, M \, \|x\|_\ga \ep^{\ga+\la-1}\, (2T)^{1-\la} =\frac12,
$$
that is
$$
\ep= \lc 2 c_{\ga,\la} \, M\, \|x\|_\ga  (2T)^{1-\la} \rc^{-1/(\ga+\la-1)}.
$$
Plugging this relation into (\ref{eq:rel-epsilon-K}), we obtain the result when
$\ep<T$.

\smallskip

Finally, $T<\ep$ if and only if $T^{\ga}<[2^{2-\la}c_{\ga+\la}M
||x||_{\ga}]^{-1}$. Thus, by inequality (\ref{numagre}), the proof is complete.

\end{proof}

\subsection{Case of a weighted delay}\label{sec:weight-delay}
In this subsection, we prove that our Hypotheses \ref{hyp:H1} and
\ref{hyp:H2} are satisfied for the weighted delay alluded to in the introduction,
that is for the function $f$ given by equation (\ref{eq:expr-f-weighted-delay}).
\begin{proposition}\label{example}
Let $\nu$ be a finite measure on $[-h,0]$ and 
$\si :\R^n\rightarrow\R^{n\times d}$ a four times differentiable 
bounded function with bounded derivatives. 
Then Hypotheses \ref{hyp:H1} and \ref{hyp:H2} are fulfilled for 
$f:\cac_1^\la([-h,0];\R^n)\to\R^{n\times d}$ defined by:
$$f(Z)=\si\left(\int_{-h}^0Z(\theta)\nu(d\theta)\right),$$
with $Z\in\cac^{\la}_1([-h,0];\R^n)$.
\end{proposition}

\begin{proof}
We first show that Hypothesis \ref{hyp:H1} holds. More specifically, the condition $|f(\zeta)|\le M$ being obvious in our case, we focus on the second condition of Hypothesis \ref{hyp:H1}. Let $Z_1,Z_2\in
\cac^{\la}_1([-h,0];\R^n)$. Then there is a constant $C>0$
such that 
\begin{eqnarray*}
\lefteqn{|f(Z_1)-f(Z_2)|}\\
&\le& C\left|\int_{-h}^0 \left(
Z_1(\theta)-Z_2(\theta)\right)\nu(d\theta)\right|
\le C\nu([-h,0])\left(\sup_{\theta\in[-h,0]}|Z_1(\theta)-Z_2(\theta)|\right).
\end{eqnarray*}
Therefore Hypothesis \ref{hyp:H1} is satisfied in this case.

Now we prove that $\UU^{(a)}$ is Fr\'echet differentiable 
in order to analyze  Hypothesis \ref{hyp:H2}. 
Since the map $Z\mapsto \int_{-h}^0Z(\cdot+\theta)\nu(d\theta)$ is 
easily shown to be
a bounded linear operator from $\cac^{\la}_{1}([a_1-h,a_2];\R^n)$
into $\cac^{\la}_{1}([a_1,a_2];\R^n)$,
we only need to show that 
$$
\si:\cac^{\la}_{\rho,a_1,a_2}(\R^n) \rightarrow\cac^{\la}_{\hro,a_1,a_2}
(\R^{n\times d}),
\quad\mbox{with}\quad
\hro\triangleq \si(\rho),
$$ 
is Fr\'echet differentiable in the directions
of $\cac^{\la}_{0,a_1,a_2}(\R^n)$, with derivative $[D\si(Z)\ell](t)=
\si'(Z(t))\ell(t)$. Towards this
end, we have to show that, taking $Z\in\cac^{\la}_{\rho,a_1,a_2}(\R^n)$ and
$\ell\in\cac^{\la}_{0,a_1,a_2}(\R^n)$, and setting
$$
q_{t}=\si(Z(t)+\ell(t))-\si(Z(t))-\si'(Z(t))\,\ell(t),
$$
then
\beq\label{eq:lim-qst}
\lim_{\|\ell\|_{\la,[a_1-h,a_2]}\rightarrow 0} \frac{\|q\|_{\la,[a_1-h,a_2]}}
{\|\ell\|_{\la,[a_1-h,a_2]}} = 0.
\eeq
In order to prove relation (\ref{eq:lim-qst}), define a 
function $b:\ou^2\to\R$ by:
$$
b(\la,\mu)=Z(s)+\la \ell(s) +\mu [Z(t)-Z(s)] +\la \mu \lc  \ell(t)-\ell(s)\rc.
$$
Observe then that $b(1,1)=Z(t)+\ell(t)$, $b(1,0)=Z(s)+\ell(s)$, $b(0,1)=Z(t)$
 and $b(0,0)=Z(s)$.
We will also set $H(\la,\mu)=\si(b(\la,\mu))$. Then
\begin{align*}
&\si(Z(t)+\ell(t))-\si(Z(t))-\si'(Z(t))\,\ell(t)  \\
&=\si(b(1,1))-\si(b(0,1))-\si'(b(0,1)) [b(1,1)-b(0,1)] 
=\frac12\iou \partial_{\la\la}^2 H(\la,1) [1-\la] \, d\la,
\end{align*}
and similarly, we have:
$$
\si(Z(s)+\ell(s))-\si(Z(s))-\si'(Z(s))\,\ell(s)
= \iou \partial_{\la\la}^2 H(\la,0) [1-\la] \, d\la.
$$
Hence, plugging these two relations in the definition of $q$, we end up with:
\begin{eqnarray*}
(\delta q)_{st}&=&\iou \lp \partial_{\la\la}^2 H(\la,1)-
\partial_{\la\la}^2 H(\la,0) \rp [1-\la] \, d\la  \\
&=& \iou \partial_{\la\la\mu}^3 H(\la,0) [1-\la] \, d\la
+\int_{\ou^2}\partial_{\la\la\mu\mu}^4 H(\la,\mu) [1-\la] [1-\mu]\, d\la d\mu.
\end{eqnarray*}
The calculation of $\partial_{\la\la\mu}^3 H(\la,0)$ and 
$\partial_{\la\la\mu\mu}^4 H(\la,\mu)$
is a matter of long and tedious computations, which are 
left to the reader. Let us just mention
that both expressions can be written as a sum of terms from which a typical
example is:
\beq\label{eq:example-deriv-along-a}
\si'''(b(\la,\mu)) \, \lc (\der Z)_{st}+\mu (\der Z)_{st} \rc \, \lc
 \ell(s)+\la (\der \ell)_{st} \rc 
\,(\der \ell)_{st}.
\eeq
These terms are obviously quadratic in $\ell$, and can be bounded uniformly
in $\la,\mu,s,t$ under the hypothesis $\si\in C_b^4$. Notice that, in order to
 bound
the term $|\ell(s)|$ in (\ref{eq:example-deriv-along-a}), we use the fact
 that $\ell$ has
a null initial condition, which means in particular that 
$|\ell(s)|\le (a_2-a_1+h)^\la \|\ell\|_{\la,[a_1-h,a_2]}$. This finishes the 
proof of (\ref{eq:lim-qst}). The continuity of $D\si(Z)$ and the existence of the constant $c_N$ introduced in Hypothesis~\ref{hyp:H2} are now a question of trivial considerations, and this ends the proof of our proposition.

\end{proof}

\begin{remark}
The proof of Frechet differentiability of $f$ was not necessary for 
the existence-uniqueness
result, which relied on some Lipschitz type condition. However, this 
stronger result turns out
to be useful for the Malliavin calculus part, and this is why we prove it 
here. Nevertheless, notice
that Theorem \ref{exiuniq} holds true for a $C_b^2$ coefficient $\si$.
\end{remark}

\subsection{Differentiability of the solution}
In this section we study the differentiability of the solution
of (\ref{deleq}) as a function of the integrator $x$, following closely the methodology of \cite{NS}. In particular, our differentiability result  will be
achieved with the help of the map 
$F:\cac^{\ga}_{0,0,T}(\R^d)\times \cac_{0,0,T}^{\la}(\R^n)\rightarrow
 \cac^{\la}_{0,0,T}(\R^n)$ 
given by
\begin{equation}\label{ope-sol}
[F(k,Z)]_t=Z_t-\cj_{0t}\lp f(\cz^{Z+{\tilde \xi}})\, d(x+k) \rp,
\quad t\in\ott
\end{equation}
where ${\tilde \xi}_t=\xi_0$ for $t\in[0,T]$, and ${\tilde \xi}_t=
\xi_t$ for $t\in[-h,0]$. Here
we recall that $\xi$ stands for an initial condition in
 $\cac_1^\la([-h,0])$. 
In this section the coefficient $f$ will satisfies the following:

\begin{hypothesis}\label{hyp:H3}
Set $\btt=(0,t)$, and recall that the map $\UU^{(\btt)}$ has been defined at
Lemma \ref{integrand}. We assume that $\UU^{(\btt)}:\cac^{\la}_{\xi,0,t}(\R^n)
\rightarrow\cac^\la([0,t];\R^{n\times d})$ 
is continuously Fr\'echet differentiable in the directions of 
$\cac_{0,0,t}^{\la}(\R^n)$, 
for some $\la\in(1/2,\ga)$. We call $\nabla\UU^{(\btt)}:\cac^{\la}_{\xi,0,t}(\R^n)
\rightarrow{\mathcal L}(\cac_{0,0,t}^{\la}(\R^n);\cac^\la_{0,0,t}(\R^{n\times d}))
$ its differential, where 
${\mathcal L}(\cac_{0,0,t}^{\la}(\R^n);\cac^\la_{0,0,t}(\R^{n\times d}))$ denotes 
the linear operators
from $\cac_{0,0,t}^{\la}(\R^n)$ into $\cac^\la_{0,0,t}(\R^{n\times d})$.
Moreover, for $s<t$ and $Z\in\cac_{0,0,T}^{\la}(\R^n)$,
$$[\nabla\UU^{(\btt)}(y)](Z)=[\nabla\UU^{(\bss)}(y)](Z)\quad
\hbox{\rm on}\quad [0,s],$$
where $y$ is the solution of equation (\ref{deleq}).
\end{hypothesis}
\begin{remarks}\label{rem:ojon}
\textbf{(1)}
Notice that 
we have shown, during the proof of Proposition \ref{example}, 
that the weighted
delay given by (\ref{eq:expr-f-weighted-delay}) also satisfies this last 
assumption.

\smallskip

\noindent
\textbf{(2)} 
If $Z\in\cac_{0,0,t}^{\la}(\R^n)$, then 
$$\|\nabla\UU^{(\btt)}(y)(Z)\|_{\la,[0,t]}\le
|\nabla\UU^{(\bt)}(y)|\|Z\|_{\la,[0,t]}.$$
Indeed, set ${\tilde Z}_s=Z_s$ for $s\in[0,t]$, and ${\tilde Z}_s=Z_t$
for $s>t$. Therefore Hypothesis \ref{hyp:H3} implies
\begin{equation*}
\|\nabla\UU^{(\btt)}(y)(Z)\|_{\la,[0,t]}\le
\|\nabla\UU^{(\bt)}(y)({\tilde Z})\|_{\la,[0,T]}
\le|\nabla\UU^{(\btt)}(y)|\|Z\|_{\la,[0,T]}
=\le|\nabla\UU^{(\btt)}(y)|\|Z\|_{\la,[0,t]},
\end{equation*}
and our claim is satisfied.
\end{remarks}

\vspace{0.2cm}

We are now ready to prove the differentiability properties for equation 
(\ref{deleq}):
\begin{lemma}\label{diff-sol} Under the Hypothesis \ref{hyp:H3}, the map 
 $F$ given by 
(\ref{ope-sol}) is continuously Fr\'echet differentiable.
\end{lemma}

\begin{proof}
Let us call respectively $D_1$ and $D_2$ the two directional derivatives. 
We first observe that,
for $k,g\in\cac^{\ga}_{0,0,T}
(\R^d)$ and $Z\in\cac_{0,0,T}^{\la}(\R^n)$, we have:
$$F(k+g,Z)-F(k,Z)+\int_0^{\cdot}\lc \UU^{(\bt)}(Z+{\tilde \xi})\rc_{s}dg_{s}=0.$$
In other words, the partial derivative $D_1F$ is defined by
$$
D_1F(k,Z)(g)=-\int_0^{\cdot}\lc\UU^{(\bt)}(Z+{\tilde \xi})
\rc_{s}dg_{s}=-\cj_{0\cdot}\lp [\UU^{(\bt)}(Z+{\tilde \xi})] \, dg \rp.
$$
We shall prove now that $D_1F$ is continuous: 
consider $k,\tilde k\in\cac^{\ga}_{0,0,T}(\R^d)$
and $Z,\tilde Z\in\cac_{0,0,T}^{\la}(\R^n)$. For notational sake, set also $\|\cdot\|_{\la}$
for $\|\cdot\|_{\la,[0,T]}$. Then, according to Lemma~\ref{ineq-int}, we obtain:
\begin{align*}
&\left\|D_1F(k,Z)(\eta)-D_1F({\tilde k},{\tilde Z})
(\eta)\right\|_{\la}
=
\left\|\cj\lp [\UU^{(\bt)}(Z+{\tilde \xi})-\UU^{(\bt)}({\tilde Z}+
{\tilde \xi}) ] \, d\eta_{s}\rp\right\|_{\la}\\
&\le\left\|\eta\right\|_{\ga}\left(\left\|\UU^{(\bt)}(Z+{\tilde \xi})-\UU^{(\bt)}
({\tilde Z}+{\tilde \xi})\right\|_{\infty}T^{\ga-\la}\right. \\
&\hspace{6cm}\left. +C_{\la+\ga}T^{\ga}
\left\|\UU^{(\bt)}(Z+{\tilde \xi})-\UU^{(\bt)}
({\tilde Z}+{\tilde \xi})\right\|_{\la}\right),
\end{align*}
which, owing to Hypothesis \ref{hyp:H3}, implies that $D_1F$ is continuous.

\vspace{0.2cm}

Concerning $D_2F$ we have, for $k\in\cac^{\ga}_{0,0,T}(\R^d)$,
$Z\in\cac_{0,0,T}^{\la}(\R^n)$ and $\tilde Z\in\cac_{0,0,T}^{\la}(\R^n)$,
and thanks to Theorem~\ref{prop-cj}:
\begin{eqnarray*}
\lefteqn{\left\| F(k,Z+{\tilde Z})-F(k,Z)-{\tilde Z}+
\cj \lp [\nabla\UU^{(\bt)}(Z+{\tilde \xi})]({\tilde Z})
\,d(x+k)\rp\right\|_{\la}}\\
&\le& \left\|x+k\right\|_{\ga}\left(\left\|\UU^{(\bt)}(Z+{\tilde Z}+{\tilde \xi})
-\UU^{(\bt)}(Z+{\tilde \xi})-[\nabla\UU^{(\bt)}(Z+{\tilde \xi})]({\tilde Z})
\right\|_{\infty}
T^{\ga-\la}\right.\\
&&\qquad\quad \left.+C_{\la+\ga}T^{\ga}
\left\|\UU^{(\bt)}(Z+{\tilde Z}+{\tilde \xi})-\UU^{(\bt)}(Z+{\tilde \xi})
-[\nabla\UU^{(\bt)}(Z+{\tilde \xi})]({\tilde Z})\right\|_{\la}
\right).
\end{eqnarray*}
Therefore, making use of Hypothesis \ref{hyp:H3}, we have that:
$$
 D_2F(k,Z)({\tilde Z})={\tilde Z}-
\int_0^{\cdot}\nabla\UU^{(\bt)}(Z+{\tilde \xi})({\tilde Z})_{s}d(x_{s}+k_{s}).
$$
The continuity of $D_2F$ can now be proven along the same lines as for $D_1F$, and the computational details are left to the reader for sake of conciseness.  The proof is now finished.

\end{proof} 

The following will be used to show that $D_2F(k,Z)$ is a linear 
homeomorphism.

\begin{lemma}\label{der-deq}
Let $w\in\cac^{\la}_{0,0,T}(\R^n)$, $y$ the
solution of (\ref{deleq}) and assume Hypotheses \ref{hyp:H1},
\ref{hyp:H2} and \ref{hyp:H3}
hold. Then the equation
\begin{equation}\label{eq-der}
Z_t=w_t+\int_0^t \left([\nabla\UU^{(\bt)}(y)](Z)\right)_{s}\,
 dx_{s},\quad 0\le t\le T,
\end{equation}
has a unique solution $Z$ in $\cac^{\la}_{0,0,T}(\R^n)$.
\end{lemma}

\begin{proof} Similarly to the proof of Theorem \ref{exiuniq}, we choose
$\ep\in(0,T)$ and
set
${\tilde \TT}_0:\cac_{0,0,\ep}^{\la}(\R^n)\rightarrow\cac_{0,0,\ep}^{\la}(\R^n)$
given by 
${\tilde \TT}_0(Z)=w+\cj_{0\cdot}([\nabla\UU^{({\mathbf{\ep}})}(y)](Z)\,dx)$. Then,
 Lemma \ref{ineq-int} and Remark \ref{rem:ojon}.(2) yield

\begin{eqnarray*}
\lefteqn{
\left\|{\tilde \TT}_0(Z)-{\tilde \TT}_0({\tilde Z})\right\|_{\la,[0,\ep]}}\\
&=& \left\| \cj\left([\nabla\UU^{({\mathbf{\ep}})}(y)](Z-\tilde Z)\,dx\right) 
\right\|_{\la,[0,\ep]}\\
&\le&\|x\|_{\la}\ep^{\ga-\la}\left(\left\|\nabla\UU^{({\mathbf{\ep}})}
(y)(Z-{\tilde Z})
\right\|_{\infty,[0,\ep]}+c_{\la+\ga}T^{\la}\left\|\nabla\UU^{({\mathbf{\ep}})}(y)
(Z-{\tilde Z})
\right\|_{\la,[0,\ep]}\right)\\
&\le&|\nabla\UU^{(\bt)}(y)|\ep^{\ga-\la}\|x\|_{\la}
\|Z-{\tilde Z}\||_{\la,[0,\ep]}(T^{\la}+c_{\la+\ga}T^{\la}).
\end{eqnarray*}
That is, for $\ep$ small enough there exists $0<C<1$ such that
$$ \left\|{\tilde \TT}_0(Z)-{\tilde \TT}_0({\tilde Z})\right\|_{\la,[0,\ep]}
\le C\|Z-{\tilde Z}\||_{\la,[0,\ep]}.$$
Hence, by standard contraction arguments, one can find a unique $Z^\ep\in\cac^\la_{0,0,\ep}(\R^n)$ such that
$$Z_t^\ep=w_t+\int_0^t \left([\nabla\UU^{({\mathbf{\ep}})}(y)](Z^\ep)\right)_{s}\,
 dx_{s},\quad 0\le t\le \ep.$$

Now we introduce 
${\tilde \TT}_\ep:\cac_{Z^\ep,\ep,2\ep}^{\la}(\R^n)
\rightarrow\cac_{Z^\ep,\ep,2\ep}^{\la}(\R^n)$
defined by
$${\tilde \TT}_\ep(Z)(t)=w_t-w_\ep+Z^\ep_\ep
+\int_\ep^t([\nabla\UU^{({\mathbf{2\ep}})}(y)](Z))_s\,dx_s,\quad t\in[\ep,2\ep].$$
Then, as in the beginning of this proof, we have
$$\left\|{\tilde \TT}_\ep(Z)-{\tilde \TT}_\ep({\tilde Z})\right\|_{\la,[0,2\ep]}
\le C\|Z-{\tilde Z}\||_{\la,[0,2\ep]}.$$
Therefore, there is a unique $Z^{2\ep}\in \cac^\la_{Z^\ep,\ep,2\ep}(\R^n)$
such that
$$
Z_t^{2\ep}=w_t+\int_0^t \left([\nabla\UU^{({\mathbf{2\ep}})}(y)](Z^{2\ep})
\right)_{s}\,
 dx_{s},\quad 0\le t\le 2\ep,$$
due to Hypothesis \ref{hyp:H3}.

Finally by induction, we can figure out a function 
$Z^{k\ep}\in \cac^\la_{Z^{(k-1)\ep},(k-1)\ep,k\ep}(\R^n)$
such that
$$
Z_t^{k\ep}=w_t+\int_0^t \left([\nabla\UU^{({\mathbf{k\ep}})}(y)](Z^{k\ep})
\right)_{s}\,
 dx_{s},\quad 0\le t\le k\ep.$$
Consequently, by Remark \ref{rem:ojon}.(2), it is not difficult to see that
$Z_t=Z_t^{k\ep}$ for $t\in[(k-1)\ep,k\ep]$ is the unique solution to equation
(\ref{eq-der}).
\end{proof}

\begin{proposition}\label{deri-solu}
Assume that Hypotheses \ref{hyp:H1} to \ref{hyp:H3} are satisfied. 
Let $y$ be the solution of
equation (\ref{deleq}). Then the map $h\mapsto y(x+h)$ is
Fr\'echet differentiable in the directions of $\cac_{0,0,T}^\ga(\R^d)$,  as  a 
$\cac^{\la}_{\xi,0,T}(\R^n)$-valued function.  Moreover, for $h,k\in
\cac_{0,0,T}^\ga(\R^d)$, we have 
\begin{eqnarray}\label{eq:der-sol}
\lc Dy(x )(k)\rc_t&=& \int_0^{t}\UU^{(\bt)}(y(x ))_{s}dk_{s}\nonumber\\
&&+\int_0^{t}\left[\nabla\UU^{(\bt)}(y(x ))(Dy(x )(k))
\right]_{s}dx_{s}.
\end{eqnarray}
In particular, $[Dy(x)](k)$ is an element of $\cac_{0,0,T}^\la(\R^n)$.
\end{proposition}

\begin{remark}\label{rem:deri-sol}  
Let us recall that equation (\ref{eq:der-sol}) has a unique solution, 
thanks to Lemma~\ref{der-deq}.
 \end{remark}

\begin{proof}[Proof of Proposition \ref{deri-solu}:]
 Like in \cite{NS}, the proof of this result is a consequence of the implicit function
theorem, and we only need to show that $D_2F(0,y(x)-{\tilde \xi})$
is a linear homeomorphism from $\cac_{0,0,T}^{\la}(\R^n)$ onto 
$\cac^{\la}_{0,0,T}(\R^n)$. Indeed, in this case we deduce that
$h\mapsto y(x)$ is Fr\'echet differentiable with
\begin{equation}\label{eq:31}
Dy(x)(k)= 
-\left(D_2F(h,y(x)-{\tilde \xi})\right)^{-1}\circ D_1F(h,y(x)-
{\tilde \xi})(k),
\end{equation}
which yields that (\ref{eq:der-sol}) holds.

\vspace{0.2cm}

Finally, notice that $D_2F(0,y(x)-{\tilde \xi})$ is bijective and continuous
 according to  Lemmas \ref{diff-sol} and \ref{der-deq}. 
Consequently the open mapping theorem implies that the application
$D_2F(0,y(x)-{\tilde \xi})$ is also a homeomorphism.

\end{proof}

\smallskip

Interestingly enough, in the particular case of the weighted delay of Section \ref{sec:weight-delay}, one can also derive a linear equation for the derivative $[Dy(x)]_t$, seen as a Hölder-continuous function.
\begin{proposition}\label{prop:maingo}Let  $\si$ and $\nu$ be as in
 Proposition \ref{example}. Let also $f$ and $y$ be defined  by (\ref{eq:expr-f-weighted-delay}) and 
(\ref{deleq}), respectively. Assume that $\nu$ is absolutely continuous 
with respect to the Lebesgue measure with Radon-Nykodim derivative in 
$L^p([-h,0])$ for $p>1/(1-\ga)$.  Then, for
$i\in\{1,\ldots,n\}$ and $k\in\cac^\la_{0,0,T}(\R^n)$, we have
$$Dy^i_t(x)(k)=\sum_{j=1}^d\int_0^t\Phi_t^{ij}(r)dk^j_r,$$
where, for $j\in\{i,\ldots,d\}$
 and $i\in\{1,\ldots,n\}$, $\Phi^{ij}$ is defined by the equation
\begin{equation}\label{eq:phim}
\Phi^{ij}_t(r)=(\UU^{(\bt)}(y))_t^{ij}+
\sum_{m=1}^n\sum_{l=1}^d\int_r^t\left(([\nabla\UU^{(\bt)}(y)]^m)^{il}
(\Phi^{mj}(s))\right)_s dx^l_s,\quad 0\le r\le t\le T,
\end{equation}
and
$\Phi_t(r)=0$ for all $0\le t<r\le T.$
\end{proposition}

\begin{remark}\label{rem:phim} Note that, for each $s\in[0,T]$ equation
(\ref{eq:phim}) has a unique solution in $\cac^\la([s,T];\R^n)$
due to Lemma \ref{der-deq}. 
\end{remark}

\begin{proof}[Proof of Proposition \ref{prop:maingo}] 
In order to avoid cumbersome matrix notations, we shall prove this result for $n=d=1$: notice that an easy consequence of the proof of Proposition \ref{example} is that in our particular case,
\begin{equation}\label{eq:evalua}
\left[\nabla
\UU^{(\bt)}(Z)(k)\right]_t=\si'\left(\int_{-h}^0Z_{t+\theta}\,\nu(d\theta)
\right)\left(\int_{-h}^0 k_{t+\theta}\,\nu(d\theta)\right).
\end{equation}
Set now $q_t=\si(\int_{-h}^0 y_{t+\theta} \,\nu(d\theta))$ and $q'_t=\si'(\int_{-h}^0 y_{t+\theta}\,\nu(d\theta))$, and write $y=y(x)$. Then equation (\ref{eq:der-sol}) can be read as:
\begin{equation}\label{eq:exp-Dy-k-onedim}
[Dy(k)]_t = \iot q_s \, dk_s+ U_t,
\quad\mbox{with}\quad
U_t=\iot q'_s \lp \int_{-h}^{0} [Dy(k)]_{s+\theta} \, \nu(d\theta) \rp \, dx_s.
\end{equation}
The Fubini type relation given at Lemma \ref{prop:25f} allows then to show, as in \cite[Proposition 4]{NS}, that 
\begin{equation}\label{eq:def-Phi}
[Dy(k)]_t=\int_0^t\Phi_t(r)dk_r,
\end{equation}
for a certain function $\Phi$, $\la$-Hölder continuous in all its variables. In order to identify the process $\Phi$, plug relation (\ref{eq:def-Phi}) into equation (\ref{eq:exp-Dy-k-onedim}) and apply Fubini's theorem, which yields
$$
U_t=\int_{-h}^{0} \nu(d\theta) 
\iot q'_s \lp \int_{0}^{(s+\theta)_{+}} \Phi_{s+\theta}(r) \, dk_r \rp \, dx_s.
$$
It should be noticed that this point is where we use the fact that $\nu(d\theta)=\mu(\theta)\, d\theta$ with $\in L^\la([-r,0])$. Indeed, in order to apply Lemma \ref{prop:25f} to $x$, $k$ and $\eta\mapsto F(\eta)=\int_{-h}^{\eta}\mu(\theta)\, d\theta$, we will assume (though this is not completely optimal) that $F$ is $\gamma$-Hölder continuous. However, a simple application of Hölder's inequality yields
$$
|F(\eta_2)-F(\eta_1)|\le c |t-s|^{(p-1)/p} \, \|\mu\|_{L^{p}([-h,0])}.
$$
It is now easily seen that the condition $(p-1)/p>\gamma$ imposes $p>1/(1-\ga)$.

\smallskip

Owing now to a (slight extension of) Lemma \ref{prop:25f}, we can write
$$
U_t=\int_{-h}^{0} \nu(d\theta) 
\int_{0}^{(t+\theta)_+} m_t(r,\theta) \, dk_r,
\quad\mbox{with}\quad
m_t(r,\theta)=\int_{r-\theta}^{t} q'_s \, \Phi_{s+\theta}(r) \, dx_s.
$$
Apply Fubini's theorem again in order to integrate with respect to $k$ in the last place: we obtain
\begin{equation*}
U_t
=\iot \lp 
\int_{-[(t-r)\wedge h]}^{0} m_t(r,\theta) \, \nu(d\theta) 
 \rp \, dk_r
=\iot \lp 
\int_{-[(t-r)\wedge h]}^{0} \nu(d\theta) \int_{r-\theta}^{t} q'_s   \, \Phi_{s+\theta}(r) \, dx_s 
 \rp \, dk_r,
\end{equation*}
and going back to (\ref{eq:exp-Dy-k-onedim}), which is valid for any $\la$-Hölder continuous function $k$, we get that $\Phi_t$ is defined on  $[0,t]$ by the equation
$$
\Phi_t(r)= q_t + 
\int_{-[(t-r)\wedge h]}^{0} \lp \int_{r-\theta}^{t} q'_s   \, \Phi_{s+\theta}(r) \, dx_s \rp 
\, \nu(d\theta),
$$
and $\Phi_t(r)=0$ if $r>t$. A last application of Fubini's theorem allows then us to recast the above equation as
$$
\Phi_t(r)= q_t + 
\int_r^t q'_s \lp \int_{-[h\wedge(s-r)]}^{0} \Phi_{s+\theta}(r) \, \nu(d\theta) \rp \, dx_s.
$$
Notice now that, if $\theta\le-(s-r)$ in the above equation, then $s+\theta\le r$, which means that $\Phi_{s+\theta}(r)=0$. Hence, we end up with an equation of the form
$$
\Phi_t(r)= q_t + 
\int_r^t q'_s \lp \int_{-h}^{0} \Phi_{s+\theta}(r) \, \nu(d\theta) \rp \, dx_s,
$$
which is easily seen to be of the form (\ref{eq:phim}).

\end{proof}

\subsection{Moments of linear equations}

In order to obtain the regularity of the density for  equation (\ref{deleq}),
 we should
bound the moments of the solution to equation (\ref{eq-der}). 
This is obtained in the following
proposition:
\begin{proposition}\label{prop:nueint-momentos-lineales}Let
${\tilde f}$ be a mapping from $\cac^{\la}_{\xi,0,T}(\R^n)$ into
the linear operators from $\cac^{\la}_{0,0,T}(\R^n)$
into $\cac^\la([0,T];\R^{n\times d})$ such that, for $0\le a<b\le T$,
${\tilde y}\in\cac^{\la}_{\xi,0,T}(\R^n)$ and ${\tilde z}\in
\cac^{\la}_{0,0,T}(\R^n)$,
\begin{itemize}
\item[(1)] $\|{\tilde f}({\tilde y}){\tilde z}\|_{\infty,
[a,b]}\le M\|{\tilde z}\|_{\infty,[a-h,b]}$.
\item[(2)]  $\|{\tilde f}({\tilde y}){\tilde z}\|_{\la,[a,b]}\le
M\|{\tilde z}\|_{\la,[a-h,b]}+ M\|{\tilde y}\|_{\la,[a-h,b]}
\|{\tilde z}\|_{\infty,[a-h,b]}.
$
\end{itemize}

 Also let 
$y$ be the solution of the equation (\ref{deleq}), $w\in\cac^{\la}_{0,0,T}
(\R^n)$
and $z\in\cac^{\la}_{0,0,T}(\R^{n})$  the solution of the equation
$$z_t=w_t+\int_0^t({\tilde f}(y)z)(t)dx_t,\quad t\in[0,T].
$$
Then
$$
\|z\|_{\la,[0,T]}\le 
c_1 \|w\|_{\la,[0,T]}D_{\ga,\la}^2 e^{c_2 D_{\ga,\la}},
$$
for two strictly positive constants $c_i=c_i(T,\ga,\la,M)$, $i=1,2$ and
$$
D_{\ga,\la}=
(\|\xi\|_{\la}\|x\|_{\ga})^{1/(\ga+\la)}+\|x\|_{\ga}^{1/\ga}+
\|x\|_{\ga}^{(2\la+\ga-1)/((\ga+\la)(\ga+\la-1))}.
$$
\end{proposition}

\begin{remarks}\label{rmk:moments-lin}
\textbf{(1)}
Observe that if $f$ is as in Proposition \ref{example} and
${\tilde f}=\nabla\UU^{(\bt)}$, then straightforward calculations show
that Conditions (1) and (2)
in the Proposition are satisfied. 

\smallskip

\noindent
\textbf{(2)}
The fact that $z_0=0$
implies that
$$ 
\|z\|_{\infty,[0,T]}\le 
c_1 T^\la\|w\|_{\la,[0,T]} D_{\ga,\la}^2 e^{c_2 D_{\ga,\la}}.
$$

\smallskip

\noindent
\textbf{(3)}
Let $\la=\ga$. Then
$(\ga+2\la-1)/((\ga+\la)(\ga+\la-1))$  in
Proposition \ref{prop:nueint-momentos-lineales} is smaller than 2 for
$\ga>H_0$, where $H_0=(7+\sqrt{17})/16\approx 0.6951$. This is the threshold above which our general delay equation will admit a smooth density.

\smallskip

\noindent
\textbf{(4)}
The unusual threshold $H_0$ above stems from the continuous dependence of the solution on its past, represented by the measure $\nu$. In case of a discrete delay of the form $\si(y_t,y_{t-r_1},\ldots,y_{t-r_q})$, we shall see that all our considerations are valid for any $H>1/2$.
\end{remarks}

\begin{proof}[Proof of Proposition \ref{prop:nueint-momentos-lineales}]
We first consider two generic positive numbers $k\in\N$ and $\ep$, such that 
$(k+1)\ep\le T$.
Then Theorem \ref{prop-cj}, point (2), and Conditions (1) and (2) imply
\begin{eqnarray*}
\lefteqn{
\|z-w\|_{\la,[k\ep,(k+1)\ep]}}\nonumber\\
&\le&\|{\tilde f}(y)z\|_{\infty,[k\ep,(k+1)\ep]}\|x\|_{\ga}\ep^{\ga-\la}
+
c_{\ga,\la}\|{\tilde f}(y)z\|_{\la,[k\ep,(k+1)\ep]}\|x\|_{\ga}\ep^{\ga}
\nonumber\\
&\le&M\|z\|_{\infty,[0,(k+1)\ep]}\|x\|_{\ga}
\ep^{\ga-\la}\nonumber\\
&& +c_{\ga,\la}M\|x\|_{\ga}
\left(\|z\|_{\la,[0,(k+1)\ep]}+\|z\|_{\infty,[0,(k+1)\ep]}
\|y\|_{\la,[0,T]}\right)\ep^{\ga}.
\end{eqnarray*}
The following (arguably non optimal) bound on $\|z\|_{\infty,[0,(k+1)\ep]}$ can now be easily verified by induction:
\begin{equation*}
\|z\|_{\infty,[0,(k+1)\ep]}\le 
\sum_{i=1}^{k+1}
2^{k+1-i}\|z-z_{(i-1)\ep}\|_{\infty,[(i-1)\ep,i\ep]}
\le 
\sum_{i=1}^{k+1}
2^{k+1-i}\|z\|_{\la,[(i-1)\ep,i\ep]}.
\end{equation*}
This yields
\begin{eqnarray}\label{ine:des-maestra}
\lefteqn{
\|z-w\|_{\la,[k\ep,(k+1)\ep]}}\nonumber\\
&\le&M\|x\|_{\ga}\ep^{\ga}\left(\sum_{i=1}^{k+1}
2^{k+1-i}\|z-z_{(i-1)\ep}\|_{\la,[(i-1)\ep,i\ep]}\right)\nonumber\\
&& 
+c_{\ga,\la}M\|x\|_{\ga}\ep^\ga\left(
\|z\|_{\la,[0,k\ep]}+\|z\|_{\la,[k\ep,(k+1)\ep]}\right)\nonumber\\
&&+c_{\ga,\la}M\|x\|_{\ga}\|y\|_{\la,[0,T]}\ep^{\ga+\la}\left(\sum_{i=1}^{k+1}
2^{k+1-i}\|z-z_{(i-1)\ep}\|_{\la,[(i-1)\ep,i\ep]}\right).
\end{eqnarray}
Now the proof can be split in three steps.

\vspace{0.2cm}

\noindent\textit{Step 1. Bounds depending on $\ep$.} Let 
\begin{equation}\label{eq:exp-epsilon}
\ep=
(T+[6M\|x\|_{\ga}(1+c_{\ga,\la})]^{1/\ga}+
[6M\|x\|_{\ga}c_{\ga,\la}\|y\|_{\la,[0,T]}]^{1/(\ga+\la})^{-1}\wedge T.
\end{equation}
Note that in this case, inequality
(\ref{ine:des-maestra}) yields
\begin{eqnarray}\label{ine:des-ep}
\lefteqn{
\|z\|_{\la,[k\ep,(k+1)\ep]}}\nonumber\\
&\le& 2\|w\|_{\la,[k\ep,(k+1)\ep]}+
M\|x\|_{\ga}
\ep^{\ga}\left(\sum_{i=1}^{k}
2^{k+2-i}\|z\|_{\la,[(i-1)\ep,i\ep]}\right)
\nonumber\\
&& +c_{\ga,\la}M\|x\|_{\ga}\ep^\ga\left(
2\|z\|_{\la,[0,k\ep]}+
\ep^{\la}\|y\|_{\la,[0,T]}\sum_{i=1}^{k}
2^{k+2-i}\|z\|_{\la,[(i-1)\ep,i\ep]}\right)\nonumber\\
&\le& 2\|w\|_{\la,[k\ep,(k+1)\ep]}\nonumber\\
&&+
\sum_{i=1}^{k}
2^{k+2-i}\|z\|_{\la,[(i-1)\ep,i\ep]}
\left(M\|x\|_{\ga}\ep^{\ga}+c_{\ga,\la}M\|x\|_{\ga}\ep^\ga+
c_{\ga,\la}M\|x\|_{\ga}\ep^{\ga+\la}\|y\|_{\la,[0,T]}
\right)\nonumber\\
&\le& 2\|w\|_{\la,[k\ep,(k+1)\ep]}+\sum_{i=1}^{k}
2^{k+1-i}\|z\|_{\la,[(i-1)\ep,i\ep]},
\end{eqnarray}
where we have used (\ref{eq:exp-epsilon}) in the last step.

\vspace{0.2cm}

\noindent\textit{Step 2. Bounds for $\|z\|_{\la,[k\ep,(k+1)\ep]}$.} 
Here we will use induction on $k$ to show that
\begin{equation}\label{cotanueva}
\|z\|_{\la,[(i-1)\ep,i\ep]}\le\sum_{j=1}^{i}2^{2i+1-2j}\|
w\|_{\la,[(j-1)\ep,j\ep]}.\end{equation}

By (\ref{ine:des-ep}) we have that this inequality holds for $i=1$.
Therefore we can assume that (\ref{cotanueva}) holds for any positive 
integer $i$ less o equal than $k$ to show that it is also true for
$i=k+1$.

The inequalities (\ref{ine:des-ep}) and (\ref{cotanueva})
lead us to write
\begin{eqnarray*}
\lefteqn{
\|z\|_{\la,[k\ep,(k+1)\ep]}}\\
&\le&2\|w\|_{\la,[k\ep,(k+1)\ep]}+\sum_{i=1}^{k}
2^{k+1-i}\sum_{j=1}^i2^{2i+1-2j}\|w\|_{\la,[(j-1)\ep,j\ep]}\\
&\le&2\|w\|_{\la,[k\ep,(k+1)\ep]}+\sum_{j=1}^{k}
\|w\|_{\la,[(j-1)\ep,j\ep]}2^{k+2-2j}\sum_{i=1}^{k}2^i\\
&\le&2\|w\|_{\la,[k\ep,(k+1)\ep]}+\sum_{j=1}^{k}
\|w\|_{\la,[(j-1)\ep,j\ep]}2^{2k+3-2j}.
\end{eqnarray*}
Now it is easy to see that (\ref{cotanueva}) also holds for $i=k+1$.

\noindent\textit{Step 3. Final bound.} 
Let $k_0$ such that $k_0\ep<T<(k_0+1)\ep$. Then, by Step 2 we have
\begin{multline*}
\|z\|_{\la,[0,T]}
\le\|w\|_{\la,[0,T]}\sum_{k=1}^{k_0}\sum_{j=1}^k2^{2k+1-2j}\\
\le\|w\|_{\la,[0,T]}(k_0)^2 2^{2k_0+1}
\le\|w\|_{\la,[0,T]}(2T/\ep)^2 2^{2T\ep^{-1}+3}.
\end{multline*}
Thus the proof is finished by plugging relation (\ref{eq:exp-epsilon}) into the last expression, and invoking Proposition 
\ref{prop:moments-solution-delay-eq}.
\end{proof}

The following result is a slight extension of Proposition \ref{prop:nueint-momentos-lineales}, allowing to take into account the case of constant but non vanishing functions.

\begin{corollary}\label{cor:nueint-momentos-lineales}Let
${\tilde f}$, $D_{\ga,\la}$, $w$ and $y$ be as in Proposition 
\ref{prop:nueint-momentos-lineales}. Furthermore, assume
that  ${\tilde f}$ is 
a mapping from $\cac^{\la}_{\xi,0,T}(\R^n)$ into
the linear operators from the constant functions on $[-h,T]$ 
into $\cac^\la([0,T];\R^{n\times d})$ satisfying the Conditions (1) and (2) of
 Proposition
\ref{prop:nueint-momentos-lineales} when ${\tilde z}$ is a constant function.
Then the solution of the equation
$$z_t=c+w_t+\int_0^t({\tilde f}(y)z)(t)dx_t,\quad t\in[0,T],$$
satisfies the inequality
$$
\|z\|_{\la,[0,T]}\le 
c_1 \left\|w + \int_0^{\cdot}({\tilde f}(y){\tilde c})(t)dx_t  
\right\|_{\la,[0,T]} D_{\ga,\la}^2 \, e^{c_2 D_{\ga,\la}},
$$
where $\tilde c$ stands for the constant function ${\tilde c}_t\equiv c$.
\end{corollary}

\begin{proof}
The proof is an immediate consequence of Proposition 
\ref{prop:nueint-momentos-lineales}. Indeed,
we only need to observe that
$$z_t-{\tilde c}_t=w_t+\int_0^t({\tilde f}(y){\tilde c})(t)dx_t
+\int_0^t({\tilde f}(y)(z-{\tilde c}))(t)dx_t,\quad t\in[0,T],$$
where ${\tilde c}(t)=c$, $t\in[0,T]$.
\end{proof}

\section{Delay equations driven by a fractional Brownian motion}\label{sec:appli-mallia}

Here we consider the Young stochastic delay equation
\begin{eqnarray}\label{stdeq}
y_t&=&\xi_0+\int_0^tf(\cz^y_t)dB_t,\quad 0\le t\le T,\nonumber\\
\cz^y_0&=&\xi,
\end{eqnarray}
where  $B=\{B_t;0\le t\le T\}$ is a $d$-dimensional fractional Brownian motion
(fBm) with 
parameter $H\in(1/2,1)$. The coefficient $f$ satisfies Hypotheses
\ref{hyp:H1}-\ref{hyp:H3}
 and $\xi$ is a given deterministic function in 
$\cac_1^{\ga}([-h,0];\R^n)$, for some $\la<\ga<H$.
 Remember that $\la\in(1/2,H)$ is introduced at the beginning of Section \ref{sec:young-delay-eq}.

The fBm $B$ 
is a centered Gaussian process with the covariance 
$$R_H(t,s)\delta_{i,j}=E(B_s^iB^j_t)=
\frac12 \delta_{i,j} (s^{2H}+t^{2H}-|t-s|^{2H}).$$
In particular, $B$ has $\nu$-H\"older continuous paths for any exponent 
$\nu<H$. Consequently, from Theorem \ref{exiuniq} and  Hypothesis
\ref{hyp:H1}-\ref{hyp:H3},
 equation (\ref{stdeq})
has a unique $\cac_{\xi,0,T}^{\la}(\R^n)$-pathwise solution.

Here, our main goal is to analyze the existence of a smooth density of the 
solution of  equation (\ref{stdeq}). This will be done via the Malliavin 
calculus or stochastic calculus of variations.
\subsection{Preliminaries on Malliavin calculus}
In this subsection we introduce the framework and the results
that we use in the 
remaining of this paper. Namely, we give some tools of the Malliavin calculus
 for fractional Brownian motion. Towards this end, we suppose that the
reader is familiar with the basic facts of stochastic analysis  
for Gaussian processes as presented, for example,
 in Nualart \cite{NL}.

Henceforth, we will consider the abstract Wiener space introduced in 
Nualart and Saussereau \cite{NS}, in order to take advantage of the relation 
between the Fr\'echet derivatives of the solution to equation (\ref{stdeq})
(see Proposition \ref{deri-solu})
and its derivatives in the Malliavin calculus sense
(see \cite{NL}, Proposition 4.1.3). This abstract Wiener space is constructed 
as follows (for a more detailed exposition of it, the reader can consult~\cite{NS}).

We assume that the underlying probability space $(\Omega,{\mathcal F},P)$ 
is such
that $\Omega$ is  the Banach space of all the continuous funtions
 $C_0([0,T];\R^d)$, which are zero at time $0$, endowed with the supremum norm.
$P$ is the only  probability measure such that the canonical process
$\{B_t;0\le t\le T\}$ is a $d$-dimensional fBm with parameter $H\in(1/2,1)$ and
the  $\sigma$-algebra ${\mathcal F}$ is the completion of the Borel 
$\sigma$-algebra of $\Omega$ with respect to $P$.

Two important tools related to the fBm $B$ are the completion  
${\HH}$ of the $\R^d$-valued step funcions ${\EE}$ with respect to the
inner product
$\langle (\1_{[0,t_1]},\ldots,\1_{[0,t_d]}),
(\1_{[0,s_1]},\dots,\1_{[0,s_d]})\rangle=
\sum_{i=1}^d R_H(s_i,t_i)$ and the isometry 
$K^*_H:{\HH}\rightarrow L^2([0,T]^d)$, which satisfies
$$
K^*_H((\1_{[0,t_1]},\ldots,\1_{[0,t_d]})=
(\1_{[0,t_1]}(\cdot)K_H(t_1,\cdot),\ldots,\1_{[0,t_d]}K_H(t_d,\cdot)),
$$
where 
$K_H(t,s)=c_Hs^{1/2-H}\int_s^t(u-s)^{H-3/2}u^{H-1/2}du$ is a kernel verifying
$$
R_H(t,s)=\int_0^{t\wedge s}K_H(t,r)K_H(s,r)dr.
$$ 
It should be noticed at this point that $K^*_H$ can be represented in the two following ways:
\begin{equation}\label{eq:rep-KH-star}
[K^*_H\vp]_t=\int_t^T \vp_r \, \partial_r K(r,t) \, dr
=c_H s^{1/2-H} [I_{T^-}^{H-1/2}(u^{H-1/2}\vp_u) ]_t,
\end{equation}
where $I_{T^-}^\al$ stands for the fractional integration of order $\al$ on $[0,T]$ (see \cite{nual-cours} for further details).

The isometry $K_H^*$
allows us to introduce the version of the Reproducing Kernel Hilbert space
$\HH_H$ 
associated with the process $B$. Namely,
Let ${\mathcal K}_H$ be given by
$$
{\mathcal K}_H:L^2([0,T];\R^d)\rightarrow \HH_H:=
{\mathcal K}_H(L^2([0,T];\R^d)), \quad
({\mathcal K}_Hh)(t)=\int_0^tK_H(t,s)h(s)ds.$$
The space $\HH$ is continuously and densely embedded in $\Omega$. Indeed,
it is not difficult to see that the operator 
$\RR_H:\HH\rightarrow \HH_H$ defined by
$$\RR_H\phi=\int_0^{\cdot}K_H(\cdot,s)(K_H^*\phi)(s)ds$$
embeds $\HH$ continuously and densely into $\Omega$, because, 
as it was pointed out in \cite{NS}, $\RR_H(\phi)$ is $H$-H\"older continuous.
Thus, we have that $(\Omega,\HH, P)$ is an abstract Wiener space.

Now we introduce the derivative in the Malliavin calculus sense 
of a random variable.
We say that a random variable $F$ is a smooth functional in ${\mathcal S}$ 
if it has the form $$F=f(B(h_1),\ldots,B(h_n)),$$
where $h_1,\ldots,h_n\in\HH$ and $f$ and all its partial derivatives have 
polynomial growth.
The derivative of this smooth fuctional is the $\HH$-valued random 
variable given by $$\DD F=\sum_{i=1}^n\frac{\partial f}{\partial x_i}
(B(h_1),\ldots,B(h_n))h_i.$$
For $p>1$, the
 operator $\DD$ is closable from $L^p(\Omega)$ into $L^p(\Omega;\HH)$
(see \cite{NL}). The closure of this operator is also denoted by 
$\DD$ and its domain by $\D^{1,p}$, which is the completion of ${\mathcal S}$
with respect to the norm
$$\|F\|_{1,p}^p=E(|F|^p)+E(\|\DD F\|_\HH^p).$$
The operator $\DD$ has the local property (i.e., $\DD F=0$ on 
$A\subset\Omega$
if $\1_AF=0$). This allows us to extend the domain of the operator $\DD$ 
as follows. We say that $F\in\D^{1,p}_{loc}$ if there is a sequence
$\{(\Omega_n, F_n), n\ge 1\}\subset{\mathcal F}\times \D^{1,p}$ such that
$\Omega_n\uparrow\Omega$ w.p.1 and $F=F_n$ on $\Omega_n$. In this case,
we define $\DD F=\DD F_n$ on $\Omega_n$.

It is known that, in the  abstract Wiener space $(\Omega,
\HH ,P)$, we can consider the   differentiability of random variable $F$
in the directions of $\HH$. That is, we say that $F$ is $\HH$-differentiable 
if for almost all $\omega\in\Omega$ and $h\in\HH$, 
the map $\ep\mapsto F(\omega+\ep\RR_Hh)$
is differentiable. The following result due to Kusuoka \cite{Ku} (see also
\cite{NL}, Proposition 4.1.3) will be fundamental in the study of the 
existence of smooth densities of the solution of equation (\ref{stdeq}).

\begin{proposition}\label{prop:kus}
 Let $F$ be an $\HH$-differentiable random variable. Then
$F$ belongs to the space $\D^{1,p}_{loc}$, for any $p>1$.
\end{proposition}

We will apply this result to the solution of equation (\ref{stdeq}) as 
follows. Note that for $\vp\in\HH$, we have the inequality
$$|(\RR_H\vp)^i(t)-(\RR_H\vp)^i(s)|=\left(E[|B^i_t-B^i_s|^2]\right)^{1/2}
\|\vp\|_{\HH}\le \|\vp\|_{\HH} |t-s|^H.$$
Consequently, Proposition \ref{deri-solu} (see also Lemma \ref{nose} below)
implies that the random
variable $y_t$ defined in  equation (\ref{stdeq}) is also $\HH$-differentiable,
which, together with Proposition~\ref{prop:kus}, yields that
$y^i_t$ belongs
to $\D^{1,p}_{loc}$ for every $t\in[0,T]$, $p>1$ and $i\in\{1,\ldots,n\}$.
Moreover, the relation between the $\HH$-derivative and  $\DD$ is given by
(see also Lemma \ref{rel-ders}),
\begin{equation}\label{eq:relders}
\langle \DD y_t^i,h\rangle_\HH=\frac{d}{d\ep}y^i_t
(\omega+\ep\RR_{H}h)|_{\ep=0},\quad h\in\HH.
\end{equation}.

More generally, if $\omega\mapsto X(\omega)$ is infinetely Fr\'echet
diferentiable in the directions of $\cac^\la_{0,0,T}(\R)$, then for a smooth random variable $X$, then
\begin{eqnarray*}
\lefteqn{
\langle \DD^n X,h_1\otimes\cdots\otimes h_n \rangle_{\HH^n}}\\
&\quad&=
D_{\RR_H h_1,\ldots,\RR_H h_n}X=
\frac{\partial}{\partial\ep_1}\ldots\frac{\partial}{\partial\ep_n}
X(\omega+\ep_1\RR_{h_1}+\ldots+\ep_n\RR_{h_n})|_{\ep_1=\ldots=\ep_n=0}.
\end{eqnarray*}
\subsection{Existence of the density of the solution}

In this section we establish that, for each
$t\in[0,T]$, the random variable $y_t$ 
introduced in equation (\ref{stdeq}) has a density.

Let us start with two important technical tools. The first one relates the derivative of the vector-valued quantity $y_t$ with the derivative of $y$ as a function.

\begin{lemma}\label{nose}
Let $y$ be the solution of (\ref{stdeq})
and $t\in[0,T]$. Then almost surely, $h\mapsto y_t(B+h)$
is Fr\'echet differentiable from $\cac^{\la}_{0,0,T}(\R^d)$
into $\R^n$. Furthermore
$$Dy_t(B)(h)=\left[Dy(B)(h)\right]_{t}.$$
\end{lemma}

\begin{proof} The proof is an immediate consequence of
\begin{eqnarray*}
\lefteqn{\left| y_t(x +h)-y_t(x )-\left(Dy(x )(h)\right)(t)\right|}\\
&=&\left| y_t(x +h)-y_t(x )-\left(Dy(x )(h)\right)(t)\right.\\
&&\quad\left.
-y_0(x +h)-y_0(x )-\left(Dy(x )(h)\right)(0)\right|\\
&\le&\left\|y(x +h)-y(x )-Dy(x )(h)\right\|_{\la} t^{\la},
\end{eqnarray*}
with $x,h\in\cac^{\la}_{0,0,T}(\R^d)$.

\end{proof}

\begin{lemma}\label{rel-ders}
Let $y$ be the solution of (\ref{stdeq}). Then $y^i_t$ belongs
to $\D^{1,2}_{loc}$ for every $t\in[0,T]$ and $i\in\{1,\ldots,n\}$. Moreover,
for $h\in \HH$, we have
\begin{equation}\label{eq:rel-ders}
\langle \DD y_t^i,h\rangle_{\HH}=\left[Dy^i(B)(\RR_{H}h)\right]_{t}.
\end{equation}
\end{lemma}

\begin{proof} By Proposition \ref{prop:kus} and Lemma \ref{nose},
 we have already shown that 
$y^i_t$ is in $\D^{1,2}_{loc}$ for every $t\in[0,T]$ and $i\in\{1,\ldots,n\}$.

Furthermore, by (\ref{eq:relders}) and Lemma 4.2, we have 
$$\langle\DD y^i_t,h\rangle_{\HH}=D_{\RR_H h} y^i_t=
Dy^i_t(B)(\RR_H h)=
\left(Dy^i(B)(\RR_H h)\right)(t).$$
Thus, the proof is complete.
\end{proof}

We now use the ideas of Nualart and
Saussereau \cite{NS} to state one of the  main results of this section:

\begin{theorem}\label{den-sol}
Let us assume that Hypotheses \ref{hyp:H1}-\ref{hyp:H3}
 hold, recall that $\xi$ is the (functional) initial condition of equation (\ref{stdeq}), and assume that the space
spanned by $\{(f(\xi)^{1j},\ldots,f(\xi)^{nj});$ 
$1\le j\le d\}$ is $\R^n$. Then for $t\in(0,T]$, the random variable $y_t$ given by
(\ref{stdeq}) is absolutely continuous with respect to the Lebesgue 
measure on $\R^n$.
\end{theorem}

\begin{proof} As in \cite{NS} (proof of Theorem 8), we have that
$y^i_t$ belongs to $\D^{1,2}_{loc}$. Therefore we only need to see
that the Malliavin covariance matrix
\begin{equation}\label{eq:def-mallia-matrix}
Q_t^{ij}:=\langle \DD y^i_t,\DD y^j_t\rangle_{\HH}
\end{equation}
is invertible almost surely.

For $v\in \R^n$, following \cite{NS} (proof of Theorem 8),
we have
$$v^TQ_tv=\sum_{m=1}^{\infty}\left|\langle Dy(B)(\RR_{H}h_m)(t),v\rangle_{\R^n}
\right|^2,$$
where $\{h_n, m\ge 1\}$ is a complete orthonormal system of $\HH$.

Now assume that the Malliavin matrix $Q_t$ 
is not almost surely invertible. Then, on the set of strictly positive probability where $Q_t$ is not invertible,
there exists $v\in\R^n$, $v\neq0$ such that $v^TQ_tv=0$. Moreover, recalling our notation (\ref{ope-sol}), it is clear from equation (\ref{eq:31}) that $D_2F(k,Z)$ is a linear homomorphism. Hence, we obtain that
\begin{eqnarray*}
0&=&\langle D_1F(0,y(B-{\tilde\xi}))(\RR_{H}h_m)(t),v_0\rangle_{\R^n}\\
&=&-\left\langle\int_0^t\UU^{(\bt)}(y(B))_{s}d\RR_{H}h_{m}(s),v_0\right\rangle_{\R^n}\\
&=& -\sum_{i=1}^n\sum_{j=1}^dv_0^i\int_0^t\left(\UU^{(\bt)}(y(B))\right)^{ij}
_{s}d\RR_{H}h^j_{m}(s)\\
&=&-\sum_{i=1}^n\langle v_0^i\left(\UU^{(\bt)}(y(B))\right)^{i}\1_{[0,t]},
h_m\rangle_{\HH},\quad \mbox{\rm for all }\ m\ge0,\end{eqnarray*}
where the last equality follows from \cite{NS}. For $t>0$, taking into account the definition of $\UU^{(\bt)}$ given at Lemma \ref{integrand}, we obtain that $\sum_{i=1}^nv_0^if^{ij}(\xi)=0,$ which contradicts the fact
that  $\R^n$ coincides with the space 
spanned by 
$$
\{(f(\xi)^{1j},\ldots,f(\xi)^{nj});\ 
1\le j\le d\}.$$ So we have that the Malliavin matrix $Q_t$ is invertible
for any $t\in(0,T]$, as we wished to prove.
\end{proof}

\subsection{Smoothness of the density of the solution} 

In order to avoid lengthy lists of hypothesis on our coefficients, we focus in this section on the example of the weighted delay treated at Section \ref{sec:weight-delay}.
As usual in the stochastic analysis context, we study the smoothness of the density of the random variable under consideration by bounding the $L^{-p}$ moments of its Malliavin matrix. Towards this aim, it will be useful to produce an equation solved by the Malliavin derivative of the solution $y_t$ of equation (\ref{stdeq}). This is contained in the following Lemma:
\begin{lemma}\label{lem:y-D-1p}
Under the conditions of Proposition  \ref{prop:maingo}, let $y$ be the solution to equation~(\ref{stdeq}). Assume furthermore that $B$ is a fBm with Hurst parameter $H>H_0$, where $H_0$ is defined at Remark \ref{rmk:moments-lin}. Then $y_t\in\D^{1,p}$ for any $p\ge 1$, and $\Phi_t(r):=\cd_r y_t$ is the unique solution to the following equation:
\begin{equation}\label{eq:dcp-Phi}
\Phi_t(r)= [\cu^{(\bt)}(y)]_t + V_t(r),
\ \mbox{where}\ 
V_t^{ij}(r)=\sum_{m=1}^n\sum_{l=1}^d\int_r^t\left(([\nabla\UU^{(\bt)}(y)]^m)^{il}
(\Phi^{mj}(s))\right)_s dB^l_s,
\end{equation}
with the additional constraint $\Phi_t(r)=0$ for all $0\le t<r\le T.$
\end{lemma}

\begin{proof}
The equation followed by $\cd y$ is a direct consequence of relation (\ref{eq:rel-ders}) and Proposition \ref{prop:maingo}. The fact that  $y_t\in\D^{1,p}$ when $H>H_0$ stems now from Proposition \ref{prop:nueint-momentos-lineales}.

\end{proof}

\smallskip

Now we are able to state the second main result of this section, for which we need an additional notation: for two a non-negative matrices $M,N\in\R^{n\times n}$, we write $M\ge N$ when the matrix $M-N$ is non-negative.
\begin{theorem}\label{theo:fmd} Let $f,\ \si$, $\nu$ and $B$ as in 
Lemma \ref{lem:y-D-1p}. Assume that $\si$ has bounded 
derivatives of any order and that 
\begin{equation}\label{eq:non-degenerate}
\si(\eta_1) \si(\eta_2)^* \ge \ep \id_{\R^n},
\quad\mbox{for all}\quad \eta_1,\eta_2\in\R^n.
\end{equation}
Then, for $t\in(0,T]$, $y_t$ has a $\cac^\infty$-density.
\end{theorem}

\begin{proof} 
The proof follows closely the lines of \cite[Theorem 3.5]{KS}, which is classical in the Malliavin calculus setting, and we shall thus proceed without giving too many details. Nevertheless, we shall divide our proof in two steps.

\smallskip

\noindent
\textit{Step 1:}
Let $Q_t$ be the Malliavin matrix of $y_t$, defined by (\ref{eq:def-mallia-matrix}). The standard conditions to verify in order to get a $\cac^\infty$ density are: (i) $y_t\in\D^\infty$,  and (ii) $[\det(Q_t)]^{-1}\in L^p$ for all $p\ge 1$. Condition (i) is obtained by iterating the derivatives of $y$, similarly to what is done in \cite{NS}, so that we will focus on point (ii).

\smallskip

In order to check that $[\det(Q_t)]^{-1}\in L^p$, we bound $P(|[\det(Q_t)]|^{-1}\ge \mu)$ for $\mu$ large enough, and invoke the fact that
$$
P\lp |[\det(Q_t)]|^{-1}\ge \mu\rp \le 
P\lp  Q_t \ngeq \frac{1}{\mu} \id_{\R^n} \rp.
$$
In the sequel of the proof, we will evaluate the right hand side of the above inequality.

\smallskip

\noindent
\textit{Step 2:}
In order to bound $Q_t$ from below, the basic idea is to use decomposition (\ref{eq:dcp-Phi}) for the Malliavin derivative of $y$. In this decomposition, the term $[\cu^{(\bt)}(y)]_t$ is bounded deterministically from below under the non-degeneracy condition (\ref{eq:non-degenerate}), while $V$ is a highly fluctuating quantity, since it is given by a stochastic integral with respect to $B$.

\smallskip

One can formalize the previous heuristic considerations in the following way:
\begin{equation*}
L_t=\left\|  \cu^{(\bt)}(y) \1_{[0,t]}\right\|_{\ch}^{2}
=\left\|  K_H^*\lp\cu^{(\bt)}(y) \1_{[0,t]}\rp\right\|_{L^2([0,t];\R^n)}^{2}.
\end{equation*}
Thanks to relation (\ref{eq:rep-KH-star}), one can show that
\begin{equation*}
L_t=c_H \sum_{l=1}^{n}
\int_0^t s^{1-2H} \int_s^t \int_s^t (r-s)^{H-3/2} (u-s)^{H-3/2} r^{H-1/2} u^{H-1/2}
\lla  q_r^* q_u, e_l \rra \, du dr ds,
\end{equation*}
where $\{e_l;\, l=1,\ldots,n\}$ stands for the canonical basis of $\R^n$, and where we have set $q_s=\si(\int_{-h}^0 y_{s+\theta} \, \nu(d\theta))$ as in the proof of Proposition \ref{prop:maingo}. Therefore, condition (\ref{eq:non-degenerate}) yields, for a constant $c$ which may change from line to line,
\begin{eqnarray*}
L_t&\ge& c\,
\ep \lp \int_0^t s^{1-2H} \int_s^t \int_s^t (r-s)^{H-3/2} (u-s)^{H-3/2} r^{H-1/2} u^{H-1/2}
 \, du dr ds \rp \id_{\R^n} \\
 &\ge& c\,
 \ep t^{2H} \id_{\R^n}.
\end{eqnarray*}

\smallskip

According to relation (\ref{eq:dcp-Phi}), it is now readily checked that
$$
Q_t\ge \frac{L_t}{2} -
\left\|  V_t\right\|_{\ch} \id_{\R^n}.
$$
Thus, for any strictly positive number $\al$, there exists a universal constant $c$ such that
\begin{equation*}
P\lp  Q_t \ngeq \frac{c \al \ep t^{2H}}{4} \id_{\R^n}\rp
\le P\lp \left\|  V_t\right\|_{\ch} \id_{\R^n} \ge  \frac{c \al \ep t^{2H}}{4} \rp
\le \lp \frac{4}{c \al \ep t^{2H}} \rp^{p} \, \frac{E\lc \left\|  V_t\right\|_{\ch}^{p} \rc}{\al^p}.
\end{equation*}
It is now enough to observe that $E[\|  V_t\|_{\ch}^{p} ]$ is a finite quantity for any $p\ge 1$, owing to Proposition \ref{prop:nueint-momentos-lineales}, to conclude the proof.

\end{proof}

\begin{remark}\label{rmk:discrete-delay}
As mentioned before, the restriction $H>H_0$ for the smoothness of the density of the random variable $y_t$ is due to the continuous dependence of our coefficient $f$ on the past of the solution. Indeed, in case of a discrete delayed coefficient of the form $\si(y_t,y_{t-r_1},\ldots,y_{t-r_q})$, with $q\ge 1$ and $r_1<\cdots<r_q\le h$, it can be seen that equation (\ref{stdeq}) can be reduced to an ordinary differential equation driven by $B$. This allows to apply the criterions given in \cite{HN}, which are valid up to $H=1/2$.

\smallskip

In order to get convinced of this fact, consider the simplest discrete delay case, that is an equation of the form 
\begin{equation}\label{eq:discrete-delay}
\xi_0+\int_0^t \si(y_t,y_{t-r})\, dB_t,\quad 0\le t\le T,
\end{equation}
with $r>0$. The initial condition of this process is given by $\xi\in\cac_1^\ga$ on $[-r,0]$, and we also assume that $\si$ and $B$ are real valued. Without loss of generality, one can assume that $T=m\, r$ for $m\in\N^{*}$. In this case, set $y(k)=\{y_{s+kr};\, s\in[0,r)\}$, and adopt the same notation for $B$. Then one can recast (\ref{eq:discrete-delay}) as
\begin{equation}\label{eq:discrete-delay2}
y_t(k)=y_{r}(k-1)+\int_{0}^{t} \si(y_u(k),y_u(k-1)) \, dB_u(k), \quad t\in[0,r], \, k\le m-1.
\end{equation}
Setting now $\mathbf{y}=(y(1),\ldots,y(m))^{t}$, $\mathbf{B}=(B(1),\ldots,B(k))^{t}$ and defining $\hat{\si} :\R^m\to\R^{m,m}$ by
\begin{equation*}
\hat{\si} (\eta(1),\ldots,\eta(m))={\rm Diag}(\si(\eta(1)),\ldots,\si(\eta(m))),
\end{equation*}
we can express (\ref{eq:discrete-delay2}) in a matrix form as
\begin{equation}\label{eq:discrete-delay3}
\mathbf{y}_t= \mathbf{y}_0 +\int_{0}^{t} \hat{\si} (\mathbf{y}_u(1),\ldots,\mathbf{y}_u(m)) \, d\mathbf{B}_u,
, \quad t\in[0,r].
\end{equation}
 This is now an ordinary equation driven by a $m$-dimensional fBm $\mathbf{B}$. Whenever $|\si(\eta)|\ge\ep>0$ and $H>1/2$, one can apply the non-degeneracy criterion of \cite{HN} in order to see that $y_t$ posesses a smooth density for any $t\in(0,T]$. The case of a vector valued original equation (\ref{eq:discrete-delay}) can also be handled through cumbersome matrix notations. As far as the case of a coefficient $\si(y_t,y_{t-r_1},\ldots,y_{t-r_q})$ is concerned, it can also be reduced to an equation of the form (\ref{eq:discrete-delay3}) by introducing all the quantities 
 \begin{equation*}
y_t(k_1,k_2,\ldots,k_r)=y_{t+\sum_{j=1}^{r}k_j(r_j-r_{j-1})},
\end{equation*}
where we have used the convention $r_0=0$.
\end{remark}

\vskip.5 cm
\noindent\textit{Acknowledgment.} Part of this work was done while 
Jorge A. Le\'on was visiting the Universit\'e Henri Poincar\'e
(Nancy) and Samy Tindel was visiting Cinvestav-IPN. Both are grateful for the hospitality of the respective institutions.


\end{document}